\documentclass[10pt]{article}

\usepackage[dvips]{graphicx}
\usepackage{makeidx}
\usepackage{amsmath}
\usepackage{amsfonts}
\usepackage[russian,english]{babel}

\newcommand {\ess}{{\mathrm{ess}}}
\newcommand {\sign}{{\mathrm{sign}}}
\newcommand {\const}{{\mathrm{const}}}
\newcommand {\dist}{{\mathrm{dist}}}

\newcommand {\etavec}{\boldsymbol{\eta}}
\newcommand {\omegavec}{{\boldsymbol{\omega}}}
\newcommand {\zetavec}{{\boldsymbol{\zeta}}}

\newcommand {\alphavec}{\boldsymbol{\alpha}}
\newcommand {\xivec}{{\boldsymbol{\xi}}}

\newcommand {\sigmavec}{{\boldsymbol{\sigma}}}
\newcommand {\epsvec}{{\boldsymbol{\epsilon}}}

\newcommand {\phivec}{\mbox{\boldmath $\phi$}}

\newcommand {\thetavec}{{\boldsymbol{\theta}}}

\newfont{\pseudocode}{cmtt10}

\newtheorem{thm}{Theorem}
\newtheorem{cor}{Corollary}
\newtheorem{lem}{Lemma}

\newtheorem{defn}{Definition}

\newtheorem{assume}{Assumption}

\newtheorem{hyp}{H}

\newcommand{\eps}{\varepsilon}

\newcommand{\s}{\mathcal{S}}

\newcommand{\bfx}{\mathbf{x}}
\newcommand{\bff}{\mathbf{f}}
\newcommand{\bfz}{\mathbf{z}}
\newcommand{\bfc}{\mathbf{c}}
\newcommand{\bfy}{\mathbf{y}}
\newcommand{\bfh}{\mathbf{h}}

\newcommand{\bfu}{\mathbf{u}}
\newcommand{\bfg}{\mathbf{g}}

\newcommand{\bfs}{\mathbf{s}}

\newcommand{\dpsi}{{\dot\psi}}
\newcommand{\pd}{{\partial}}

\newcommand{\Real}{\mathbb{R}}
\newcommand{\Natural}{\mathbb{N}}

\makeindex
\begin{document}
% for Improved Performance!

%\twocolumn

\baselineskip 0.7cm

\title{\bf Adaptation and Parameter Estimation in Systems with Unstable Target Dynamics and Nonlinear Parametrization}
%\title{\bf Adaptive Algorithms in a Finite Form for Nonlinear Dynamic Plants}
\author{Ivan Tyukin\thanks{Laboratory for Perceptual Dynamics, Brain Science Institute, RIKEN (Institute for Physical and Chemical Research)
                           , 2-1, Hirosawa, Wako-shi, Saitama, 351-0198, Japan, e-mail:
                           \{tyukinivan\}@brain.riken.jp}
%                           , Prashanth Alluvada\thanks{Laboratory for Perceptual Dynamics, RIKEN (Institute for Physical and Chemical Research)
%                           Brain Science Institute, 2-1, Hirosawa, Wako-shi, Saitama, 351-0198, Japan, e-mail:
%                           \{aprashanth\}@brain.riken.jp}
                           , Danil
                           Prokhorov\thanks{Ford Research Laboratory, Dearborn, MI,
48121, USA, e-mail: dprokhor@ford.com},  Cees van
Leeuwen\thanks{Laboratory for Perceptual Dynamics, Brain Science
Institute, RIKEN (Institute for Physical and Chemical Research)
                          , 2-1, Hirosawa, Wako-shi, Saitama, 351-0198, Japan, e-mail:
                           \{ceesvl\}@brain.riken.jp}
}
\date{\today}
\date{}
\maketitle{}
%\thanks{Laboratory for Perceptual Dynamics, RIKEN (Institute for Physical and Chemical Research)
%                           Brain Science Institute, 2-1, Hirosawa, Wako-shi, Saitama, 351-0198, Japan, e-mail: \{tyukinivan,ceesvl\}@brain.riken.jp}
%\thanks{Ford Research Laboratory, Dearborn, MI, 48121, USA, e-mail:
%dprokhor@ford.com}

\begin{abstract}
We propose a technique for the design and analysis of adaptation
algorithms in dynamical systems. The technique applies both to
systems with conventional Lyapunov-stable target dynamics and to
ones of which the desired dynamics around the target set is
nonequilibrium and in general unstable in the Lyapunov sense.
Mathematical models of uncertainties are allowed to be nonlinearly
parametrized, smooth, and monotonic functions of linear
functionals of the parameters. We illustrate with applications how
the proposed method leads to control algorithms. In particular we
show that the mere existence of nonlinear operator gains for the
desired dynamics guarantees that system solutions are bounded,
reach a neighborhood of the target set, and mismatches between the
modeled uncertainties and uncertainty compensator vanish with
time. The proposed class of algorithms can also serve as parameter
identification procedures. In particular, standard persistent
excitation suffices to ensure exponential convergence of the
estimated to the actual values of the parameters. When a weak,
nonlinear version of the persistent excitation condition is
satisfied, convergence is asymptotic. The approach extends to a
broader class of parameterizations where the monotonicity
restriction holds only locally. In this case excitation with
oscillations of sufficiently high frequency ensure convergence.

%We propose a novel technique for the design and analysis of
%algorithms for adaptive control and parameter estimation in
%dynamical systems. The technique applies to both systems with
%conventional Lyapunov-stable target dynamics and those of which
%the desired dynamics around the target set is non-equilibrium and
%in general is unstable in Lyapunov sense. Mathematical models of
%uncertainties are allowed to be nonlinearly parametrized, smooth,
%and monotonic functions of linear functionals of the parameters.
%For this class of systems our method leads to novel  control
%particular we show that merely existence of nonlinear operator
%gains for the desired dynamics guarantees that system solutions
%are bounded, reach a neighborhood of the target set, and
%mismatches between the modeled uncertainties and uncertainty
%compensator vanish with time.
%
%We  also show that these new algorithms can in principle serve as
%parameter identification procedures. In particular, standard
%persistent excitation suffices to ensure exponential convergence
%of the estimated to the actual values of the parameters. When new,
%much weaker nonlinear version of the persistent excitation
%condition is satisfied we prove that convergence is {\it
%asymptotic}. Last but not least we provide extensions of these
%algorithms to broader class of parameterizations where the
%monotonicity restriction holds only locally. In this case
%excitation with sufficiently  high-frequency of oscillations
%ensures convergence.  A practically relevant example is given in
%order to illustrate the effectiveness of the approach.
\end{abstract}

%\vspace{10mm}


 \vspace{10mm}
{ \small
{\bf Corresponding author:}{  \\Ivan Tyukin \\
                              Laboratory for Perceptual Dynamics,\\
                              RIKEN Brain Science Institute,\\
                              2-1, Hirosawa, Wako-shi, Saitama,\\
                              351-0198, Japan\\
                              phone: +81-48-462-1111 extension 7436\\
                              fax:   +81-48-467-7236\\
                              e-mail: tyukinivan@brain.riken.jp}
}
\section{Introduction}

Results in adaptive control theory and systems identification are
most frequently used in control engineering, but have potentially
a much wider significance. In particular these theories are of
great potential relevance for sciences such as physics and biology
\cite{Sontag_2004}. On the other hand, it is in these areas that
the current limitations of control theory are most strongly felt.
Whereas effective procedures are available in case the system is
static
\cite{Box},\cite{Hansen92},\cite{Wilde},\cite{Kirkpatrick83},
adequate solutions for dynamical systems have been proposed under
conditions that may not be adequate for most scientific
applications. These conditions require that systems are linear in
their parameters, the target dynamics is stable in the Lyapunov
sense, and a Lyapunov function of the target dynamics can be given
\cite{Sastry89}, \cite{Narendra89}, \cite{Kokotovich95},
\cite{Ljung_99},\cite{Eykhoff}, \cite{Bastin92},\cite{Fradkov99}.
Each of these restrictions alone is limiting the role of control
theory in the scientific arena; together they constitute the
"standard" approach that confines control theory to a limited
role, even within the realm of engineering.

Whereas in artificial system design, nonlinear parametrizations
could often be avoided, physical and biological models often
require the inclusion of nonlinearly parametrized uncertainties
\cite{Armstrong_1993},\cite{Pacejka91},\cite{Boskovic_1995},\cite{Abbott_2001},
\cite{Kitching_2000}.  Proposed solutions to the nonlinear
parametrization problem have cemented the standard approach, in
that they eliminate any hopes of escaping from the stable target
dynamics requirement. Nonlinearity is traditionally solved by
invoking dominance of the nonlinear terms \cite{Lin_2002_smooth},
\cite{Lin}. Dominance inevitably overcompensates the nonlinearity
inherent to the system. This is undesirable if the system's target
motions require such nonlinearities. It is in particularly
unhelpful, in case the system equations embody certain physical
laws or other regularities that necessitate nonlinear
parametrization. Terms overcompensated by dominance are likely to
be exactly the ones postulated by these laws and regularities. In
order to enable nonstable and in a sense more delicate target
dynamics, more gentle control is needed: one that enables a system
to reach the desired dynamical state by modification of, rather
than destroying, its intrinsic motions.

Alternatives to dominance are available, but they face a variety
of restrictions that make them appear less satisfactory. Often
they apply to a narrowly defined class of parameterizations, e.g.
convex functionals as in \cite{Fradkov79}. When a broader class of
nonlinearities is considered, for instance Hammerstain (Wiener)
models,
\cite{Narendra_66},\cite{Pawlak_91},\cite{Garulli_2002},\cite{Bai_2003},
the functions are restricted to static input (output)
nonlinearities. Perhaps the most promising approach so far
involves local linear (nonlinear) modelling techniques
\cite{Johansen_95}, \cite{Verdult_2002}, \cite{Enqvist_CDC2002}.
The resulting models, on the other hand, are not always physically
plausible. In case fairly general nonlinear state dependent
functions are allowed \cite{Cao_2003}, the class of dynamical
systems is limited to those modeled by the first-order ordinary
differential equations with nonlinear parameterized terms that are
Lipschitz in time. The last restriction but not least is that the
majority of these methods rely on the assumption of stable target
dynamics.

There are physical and biological systems, however, which do not
meet the requirement of stable target dynamics
\cite{Gauthier_1996}, \cite{Nature:Ditchburn:1952}. Multistability
and coexistence of multiple attractors
\cite{Arecchi:82},\cite{Weigel_98},\cite{Chizh_2000},\cite{Maldonado_1990}
are well-known examples where a system could, at best, be only
locally stable. In biophysics amplification by oscillatory
instability is believed to be a general mechanism of signal
detection in sensory systems \cite{PNAS:Camalet:2000}.
Furthermore, as demonstrated in \cite{Raffone_2003}, instability
(intermittency) offers a solution to the longstanding binding
problem in the biology of vision
\cite{Neuron:von_der_Malsburg:1999}. In fact, also in artificial
systems unstable target dynamics are sometimes required
\cite{Raffone_2003}, \cite{Conf:NOLCOS-04}. For instance, in
\cite{Neural_Computation:Suimetsu:2004}, the effectiveness of
using chaotic dynamics in solving the path finding problems in
robotics is shown. In computer science unstable (intermittent)
synchrony was shown to be an effective paradigm for solving the
image segmentation problems \cite{van_Leeuwen}.

Control-theoretic motivation and successful solutions to the
problem of adaptive regulation to unstable dynamical states are
provided in \cite{PRE:Moreau:2003},\cite{SysConLett:Moreau:2003}
for linear systems with linear parameterization. For nonlinear
systems with nonlinear parameterization, and, possibly, unstable
target dynamics the problems of adaptive regulation and parameter
estimation need further development.

In our present paper we aim to provide a unified tool capable of
solving the problems of adaptation and parameter estimation

$\bullet$ in the presence of nonlinear {\it state-dependent}
parameterization;

$\bullet$ with non-trivial target sets (namely, surfaces in the
systems's state space)

$\bullet$ with potentially unstable target dynamics, and therefore

$\bullet$ without requiring for knowledge (or existence) of
Lyapunov functions of the desired motions\footnote{This problem,
as mentioned in \cite{Fradkov99}, was long remained an open
theoretical challenge.}

Previous efforts to deal with nonlinearity by adopting domination
functions \cite{Lin},\cite{Lin_2002_smooth}, or low-order
mathematical models \cite{Cao_2003} have tried to address the most
general case. In contrast to this, we restrict ourselves in
advance to a certain  class of nonlinearities. This class,
however, is wide enough to cover a variety of relevant models in
physics, mechanics, physiology and neural computation. In
particular, it includes models of stiction, slip and surface
dependent friction, nonlinearities in dampers, smooth saturation,
dead-zones in mechanical systems, and nonlinearities in models of
bio-reactors \cite{Armstrong_1993}, \cite{Pacejka91},
\cite{Boskovic_1995}, \cite{Abbott_2001}, \cite{Kitching_2000}.

%When analyzing available approaches to deal with nonlinear in
%parameter state-dependent uncertainties we argue that use of the
%domination functions \cite{Lin},\cite{Lin_2002_smooth}, or
%low-order mathematical models \cite{Cao_2003}, are the consequence
%of a strategy to deal with nonlinear parametrization in the most
%general case. Therefore, we restrict these nonlinearities  to a
%certain {\it practically relevant} class, which is wide enough to
%include a variety of models in physics, mechanics, physiology and
%neural computation. In particular, it covers effects of stiction,
%slip and surface dependent friction, nonlinearities in dampers,
%smooth saturation, and dead-zones in mechanical systems,
%nonlinearities in models of bio-reactors \cite{Armstrong_1993},
%\cite{Pacejka91}, \cite{Boskovic_1995}, \cite{Abbott_2001},
%\cite{Kitching_2000}.

In order to deal with unstable and non-equilibrium target
dynamics, without invoking the knowledge (or even existence) of
the corresponding Lyapunov funcion, we employ {\it operator}
formalism in the {\it functional spaces} rather than conventional
tools\footnote{We refer here to common practice to fit the
derivative of  the goal functionals (Lyapunov candidates) to
specific algebraic inequalities leading to the property of
Lyapunov stability.}. In particular, we consider desired dynamics
in terms of input-output mappings in the specific functional
spaces. The only requirement we impose on these mappings is
existence of nonlinear operator gains that bound functional norms
of the outputs, given norm-bounded inputs. The inputs for these
mappings are the mismatches between the modeled uncertainty and a
compensator. The outputs are the state $\bfx$ and a function
$\psi(\bfx,t)$, not necessarily definite in state, which is
considered a measure of deviation from the target set. This
system-theoretic point of view allows us to formulate the problem
of adaptation as a problem of regulation of the mismatches to
specific functional spaces followed, if possible, by minimization
of their functional norm.

We show that the solution to this problem for the given class of
parameterizations does not require continuity of the corresponding
operator gains. This, in turn, suggests that stability of the
desired dynamics, which in many cases is synonymous to continuity
of the input-output operators \cite{Zames:66},  is not a necessary
requirement for our approach. Furthermore, given that
$\psi(\bfx,t)$ may not be definite, this new point of view on the
adaptation allows us to lift conventional state-space metric
restrictions on the goal functionals\footnote{Which are usually
defined as positive-definite and radially unbounded functions of
state \cite{Jac},\cite{Narendra89},\cite{Kokotovich95}}.

Under standard and intuitively clear additional hypotheses (i.e.
persistent excitation of a certain functional of state), we show
that the proposed adaptation procedures solve the problem of
parameter estimation for nonlinearly parameterized uncertainties.
In this case convergence is exponential. The estimates of the
convergence rates are based on the results of \cite{Lorea_2002}
and provided here for consistency. In case the conditions
specifying the class of nonlinear in parameter uncertainties hold
only locally, we show that sufficiently high frequency of
excitation still ensures convergence. For cases where the standard
persistent excitation property does not hold, we formulate a new
version of {\it nonlinear persistent excitation} condition
\cite{Cao_2003}. With this new property it is still possible to
show {\it asymptotic} convergence of the estimates to the actual
values of unknown parameters. Whether the convergence is
exponential is not answered in this paper.

The paper is organized as follows. Section 2 describes notations
and conventions we are using in the paper; in Section 3 we
formulate the problem. For the sake of compact exposition of our
results we restrict ourselves to systems that are affine in
control, although some non-affine cases are discussed towards the
end of the paper. Section 4 contains the main results of the
paper. We discuss several immediate extensions of the present
results in Section 5. In Section 6 we provide a practically
relevant application of our method, and Section 7 concludes the
paper.

\section{Notation}

According to the standard convention, symbol $\Real$ defines the
field of real numbers and  $\Real_{\geq c}=\{x\in\Real|x\geq c\}$,
$\Real_{+}=\Real_{\geq 0}$; symbol $\Natural$ defines the set of
natural numbers; symbol $\Real^n$ stands for a linear space
$\mathcal{L}(\Real)$ over the field of reals with
$\mathrm{dim}\{\mathcal{L}(\Real)\}=n$; $\|\bfx\|$ denotes the
Euclidian norm of $\bfx\in\Real^n$; $\mathcal{C}^k$ denotes the
space of functions that are  at least $k$ times differentiable.
Symbol $\mathcal{K}$ denotes the space of all functions $\kappa:
\Real_+\rightarrow \Real_+$ such that $\kappa(0)=0$, and that
$x'>x''$, $x',x''\in\Real_+$ implies that
$\kappa(x')-\kappa(x'')>0$. By symbol ${L}_{p}^n[t_0,T]$, where
$T>0$, $p\geq 1$ we denote the space of all functions
$\bff:\Real_+\rightarrow\Real^n$ such that
\[
\|\bff\|_{p,[t_0,T]}=\left(\int_{0}^T\|\bff(\tau)\|^{p}d\tau\right)^{1/p}<\infty
\]
Symbol $\|\bff\|_{p,[t_0,T]}$ denotes the ${L}_{p}^n[t_0,T]$-norm
of vector-function $\bff(t)$. By ${L}^n_\infty[t_0,T]$ we denote
the space of all functions $\bff:\Real_+\rightarrow\Real^n$ such
that
\[
\|\bff\|_{\infty,[t_0,T]}=\ess \sup\{\|\bff(t)\|,t \in
[t_0,T]\}<\infty,
\]
and $\|\bff\|_{\infty,[t_0,T]}$ stands for the
${L}^n_\infty[t_0,T]$ norm of $\bff(t)$.

Let $\bff: \Real^{n}\rightarrow \Real^m$ be given. Function
$\bff(\bfx): \Real^{n}\rightarrow \Real^m$ is said to be locally
bounded if for any $\|\bfx\|<\delta$ there exists constant
$D(\delta)>0$ such that the following holds: $\|\bff(\bfx)\|\leq
D(\delta)$.

Let $\Gamma$ be an $n\times n$ square matrix, then $\Gamma>0$
denotes a positive definite (symmetric) matrix, and $\Gamma^{-1}$
is the inverse of $\Gamma$. By $\Gamma\geq 0$ we denote a positive
semi-definite matrix. Symbols $\lambda_{\min}(\Gamma)$,
$\lambda_{\max}(\Gamma)$ stand for the minimal and maximal
eigenvalues of $\Gamma$ respectively. By symbol $I$ we denote the
identity matrix. We reserve symbol $\|\bfx\|_{\Gamma}^2$ to denote
the following quadratic form: $\bfx^{T}\Gamma\bfx$,
$\bfx\in\Real^n$. Notation $|\cdot|$ stands for the module of a
scalar. The solution of a system of differential equations
$\dot{\bfx}=\bff(\bfx,t,\thetavec,\bfu), \ \bfx(t_0)=\bfx_0$,
$\bfu:\Real_+\rightarrow\Real^m$, $\thetavec\in\Real^d$ for $t\geq
t_0$ will be denoted as $\bfx(t,\bfx_0,t_0,\thetavec,\bfu)$, or
simply as $\bfx(t)$ if it is clear from the context what  the
values of $\bfx_0,\thetavec$ are and how the function $\bfu(t)$ is
defined.

Let $\bfu:\Real^n\times\Real^d\times\Real_+\rightarrow\Real^m$ be
a function of state $\bfx$, parameters $\hat{\thetavec}$, and time
$t$. Let in addition both $\bfx$ and $\hat{\thetavec}$ be
functions of $t$. Then in case the arguments of $\bfu$ are clearly
defined by the context,  we will simply write $\bfu(t)$ instead of
$\bfu(\bfx(t),\hat{\thetavec}(t),t)$.

The (forward complete) system
$\dot{\bfx}=\bff(\bfx,t,\thetavec,\bfu(t))$, is said to have an
$L_{p}^m [t_0,T]\mapsto L_{q}^n[t_0,T]$, gain ($T\geq t_0$,
$p,q\in\Real_{\geq 1}\cup\infty$) with respect to its input
$\bfu(t)$ if and only if $\bfx(t,\bfx_0,t_0,\thetavec,\bfu(t))\in
L_{q}^n [t_0,T]$ for any $\bfu(t)\in L_{p}^m [t_0,T]$ and there
exists a function
$\gamma_{q,p}:\Real^n\times\Real^d\times\Real_+\rightarrow\Real_+$
such that the following inequality holds:
\[
\|\bfx(t)\|_{q,[t_0,T]}\leq
\gamma_{q,p}(\bfx_0,\thetavec,\|\bfu(t)\|_{p,[t_0,T]})
\]
Function $\gamma_{q,p}(\bfx_0,\thetavec,\|\bfu(t)\|_{p,[t_0,T]})$
is assumed to be non-decreasing  in $\|\bfu(t)\|_{p,[t_0,T]}$, and
locally bounded in its arguments.

For notational convenience when dealing with vector fields and
partial derivatives we will use the following extended notion of
Lie derivative of a function. Let it be the case that
$\bfx\in\Real^n$ and $\bfx$ can be partitioned as follows
$\bfx=\bfx_1\oplus\bfx_2$, where $\bfx_1\in\Real^q$,
$\bfx_1=(x_{11},\dots,x_{1q})^T$, $\bfx_2\in\Real^p$,
$\bfx_2=(x_{21},\dots,x_{2p})^T$, $q+p=n$, and $\oplus$ denotes
concatenation of two vectors. Define
$\bff:\Real^{n}\rightarrow\Real^n$ such that
$\bff(\bfx)=\bff_1(\bfx)\oplus\bff_2(\bfx)$, where
$\bff_1:\Real^n\rightarrow\Real^q$,
$\bff_1(\cdot)=(f_{11}(\cdot),\dots,f_{1q}(\cdot))^T$,
$\bff_2:\Real^n\rightarrow\Real^p$,
$\bff_2(\cdot)=(f_{21}(\cdot),\dots,f_{2p}(\cdot))^T$. Then symbol
$L_{\bff_i(\bfx)}\psi(\bfx,t)$, $i\in\{1,2\}$ denotes the Lie
derivative of function $\psi(\bfx,t)$ with respect to vector field
$\bff_i(\bfx,\thetavec)$:
\[
L_{\bff_i(\bfx)}\psi(\bfx,t)=\sum_{j}^{\dim{\bfx_i}}\frac{\pd
\psi(\bfx,t) }{\pd x_{ij}}f_{ij}(\bfx,\thetavec)
\]
Symbol $\sign(\cdot)$ denotes the signum-function:
\[
\sign(s)=\left\{
          \begin{array}{ll}
          1, & s >0\\
          0, & s =0\\
          -1, & s <0
          \end{array}
         \right.
\]

\section{Problem Formulation}

Let the following system be given:
\begin{eqnarray}\label{system1}
\dot{\bfx}_1&=&\bff_1(\bfx)+\bfg_1(\bfx)u,\nonumber \\
\dot{\bfx}_2&=&\bff_2(\bfx,\thetavec)+\bfg_2(\bfx)u,
\end{eqnarray}
where
\[
\bfx_1=(x_{11},\dots,x_{1 q})^T\in \Real^q
\]
\[
\bfx_2=(x_{21},\dots,x_{2 p})^T\in \Real^p
\]
\[
\bfx=(x_{11},\dots,x_{1 q},x_{21},\dots,x_{2 p})^T\in \Real^{n}
\]
$\thetavec\in \Omega_\theta\in \Real^d$ is a vector of unknown
parameters, and $\Omega_\theta$ is a closed bounded subset of
$\Real^d$; $u\in\Real$ is the control input, and functions
$\bff_1:\Real^{n}\rightarrow \Real^{q}$,
$\bff_2:\Real^{n}\times\Real^d\rightarrow \Real^{p}$,
$\bfg_1:\Real^{n}\rightarrow \Real^q$,
$\bfg_2:\Real^{n}\rightarrow\Real^{p}$ are locally bounded. Vector
$\bfx\in\Real^n$ is a state vector, and vectors $\bfx_1$, $\bfx_2$
are referred to as {\it uncertainty-independent} and {\it
uncertainty-dependent} partitions of $\bfx$, respectively.

For the sake of compactness we introduce the following alternative
description for (\ref{system1}):
\begin{eqnarray}\label{system}
\dot{\bfx}=\bff(\bfx,\thetavec)+\bfg(\bfx)u,
\end{eqnarray}
where
\[
\bfg(\bfx)=(g_{11}(\bfx),\dots,g_{1q}(\bfx),g_{21}(\bfx),\dots,g_{2
p}(\bfx))^{T}
\]
\[
\bff(\bfx)=(f_{11}(\bfx),\dots,f_{1q}(\bfx),f_{21}(\bfx,\thetavec),\dots,f_{2
p}(\bfx,\thetavec))^{T}
\]
Our goal is to derive both the control function $u(\bfx,t)$ and
estimator $\hat{\thetavec}(t)$, such that all trajectories of the
system are bounded and state $\bfx(t)$ converges to the desired
domain in $\Real^n$. In addition, we would like to find conditions
ensuring that the estimate $\hat{\thetavec}(t)$ converges to
unknown $\thetavec\in \Omega_\theta$ asymptotically. In order to
ensure boundedness of the trajectories, we should design an input
$u(\bfx,t)$ that restricts all possible motions of system
(\ref{system}) to an admissible bounded domain
$\Omega\subset\Real^n$ in the system state space, and if possible
steers trajectories $\bfx(t)$ to the specific set
$\Omega_0\subseteq\Omega$.

As a measure of closeness of trajectories $\bfx(t)$ to the desired
state we introduce the error function $\psi:\Real^n\times
\Real_+\rightarrow \Real, \ \psi\in \mathcal{C}^1$ such that
\begin{equation}\label{eq:target_set}
\Omega_0=\{\bfx(t)\in\Real^n|\psi(\bfx(t),t)=0\}
\end{equation}
In conventional theories it is usually required that function
$\psi(\bfx,t)$ satisfies  ({\it algebraic}) metric restrictions:
\begin{equation}\label{eq:standard_goal}
{\nu}_1(\|\bfx-\xivec(t)\|)\leq \psi(\bfx,t)\leq
{\nu}_2(\|\bfx-\xivec(t)\|), \ {\nu}_1, {\nu}_2\in
\mathcal{K}_\infty,
\end{equation}
where function $\xivec:\Real_+\rightarrow\Real^n$, $\xivec\in
\mathcal{C}^0$ is, for instance, the reference trajectory.
Function $\psi(\bfx,t)$ in this case serves as the Lyapunov
candidate for the controlled system under the assumption that
$\thetavec$ is known. The problem, however, is that finding such a
Lyapunov candidate is not a trivial task. Furthermore, the desired
trajectories $\xivec(t)$ as functions of time may only be
partially specified. In case no reference function is available
(e.g $\xivec(t)=0$) and the task is to steer the state $\bfx$ of
system (\ref{system1}) to a non-trivial set
$\Omega_0\subset\Real^n$, it is often difficult to find a goal
functional $\psi(\bfx,t)$ such that both (\ref{eq:target_set}) and
(\ref{eq:standard_goal}) are satisfied.
%Therefore, in our current
%study we would like to refrain from the conventional definition of
%a goal functional $\psi(\bfx,t)$ given by inequality
%(\ref{eq:standard_goal}).
In addition, in physical and nonlinear
systems the desired dynamics (e.g. dynamics of convergence of
trajectories $\bfx(t)$ to the reference $\xivec(t)$) could be
unstable in Lyapunov sense \cite{Gauthier_1996,Raffone_2003},
although it may posess certain degree of attraction
\cite{Milnor_1985}, and  bounded deviation from the reference.

In order to tackle these complex, but still possible, phenomena we
propose to replace the conventional goal functionals
(\ref{eq:standard_goal}) with new and less restrictive ones. In
particular, we propose to replace the standard norms $\|\cdot\|$
in $\Real^n$ in (\ref{eq:standard_goal}) with {\it functional
norms} $\|\bfx(t)\|_{p,[t_0,T]}$, $T\geq t_0$ in the {\it
functional spaces} $L_{p}[t_0,T]$, $T\geq t_0$, $p\in\Real_{\geq
1}\cup \infty$. In the other words, we replace algebraic
inequality (\ref{eq:standard_goal}) with {\it operator} relations.
This will allow us to keep function $\psi(\bfx,t)$ as a measure of
closeness of trajectories $\bfx(t)$ to the desired set $\Omega_0$
without imposing state-metric restrictions
(\ref{eq:standard_goal}) on the function $\psi(\bfx,t)$. On the
other hand we will be able to derive bounds for $\bfx(t)$ from the
values of functional $L_{p}^1[t_0,T]$-norms of the function
$\psi(\bfx(t),t)$. Let us formally introduce this requirement as
follows:
\begin{assume}[Target operator]\label{assume:psi}  For the given function $\psi(\bfx,t)\in \mathcal{C}^1$ the following  property holds:
\begin{equation}\label{eq:assume_psi}
\|\bfx(t)\|_{\infty,[t_0,T]}\leq
\tilde{\gamma}\left(\bfx_0,\thetavec,\|\psi(\bfx(t),t)\|_{\infty,[t_0,T]}\right)
\end{equation}
where
$\tilde{\gamma}\left(\bfx_0,\thetavec,\|\psi(\bfx(t),t)\|_{\infty,[t_0,T]}\right)$
is a locally bounded and non-negative function of its arguments.
\end{assume}
Assumption \ref{assume:psi} can be interpreted as a sort of {\it
unboundedness observability} property \cite{Jiang_1994} of system
(\ref{system1}) with respect to the ``output" function
$\psi(\bfx,t)$. It can also be viewed as a {\it bounded input -
bounded state} assumption for system (\ref{system1}) along the
constraint
$\psi(\bfx(t,\bfx_0,t_0,\thetavec,u(\bfx(t),t)),t)=\upsilon(t)$,
where signal $\upsilon(t)$ serves as the new input\footnote{If,
however, boundedness of the state is not required  explicitly
(i.e. it is guaranteed by additional control or follows from the
physical properties of the system itself), Assumption
\ref{assume:psi} can be removed from the statements of our
results.}. In order to illustrate this consider the equations of a
spring-mass system with nonlinear damping:
\begin{equation}\label{eq:system_spring}
\begin{split}
\dot{x}_1&=x_2 \\
\dot{x}_2&=k_0 x_1 + f(x_2,t)+u(t), \ k_0<0
\end{split}
\end{equation}
where $f:\Real\times\Real_+\rightarrow\Real$, $f(\cdot,\cdot)\in
\mathcal{C}^1$ is the nonlinear time-varying damping term.
Equations of this type arise in broad areas of engineering ranging
from  active suspension control \cite{Chantra_1999} to
 haptic interfaces \cite{Okamura2002} and identification
of the muscle dynamics \cite{Wu_1990}. Let the desired dynamics of
(\ref{eq:system_spring}) be an exponentially fast convergence of
$x_1(t)$, $x_2(t)$ to the origin. This requirement is satisfied
for the following target set:
\[
\Omega_0=\{\bfx\in\Real^2|\bfx: x_1+ \lambda x_2=0 \},  \
\lambda\in\Real, \lambda>0
\]
Therefore, function $\psi(\bfx,t)$ could be chosen as
$\psi(\bfx,t)=x_1+\lambda x_2$. Let $\psi(\bfx(t),t)\in
L_\infty^1[t_0,T]$, i. e. $x_1(t) + \lambda x_2(t)=\upsilon(t)$,
$\upsilon(t)\in L_\infty^1[t_0,T]$. An  equivalent description of
system (\ref{eq:system_spring}) in accordance with this constraint
is given by
\begin{equation}\label{eq:BIBO_spring}
\dot{x}_1=-\lambda^{-1}x_1 + \lambda^{-1}\upsilon(t), \ \lambda
x_2(t)+x_1(t)=\upsilon(t)
\end{equation}
It is clear that system (\ref{eq:BIBO_spring}) has the bounded
input - bounded state property with respect to input $\upsilon(t)$
as $\|x_1(t)\|_{\infty,[t_0,T]}\leq |x_1(t_0)| +
\|\upsilon(t)\|_{\infty,[t_0,T]}$ and
$\|x_2(t)\|_{\infty,[t_0,T]}\leq\lambda^{-1}(
\|x_1(t)\|_{\infty,[t_0,T]}+\|\upsilon(t)\|_{\infty,[t_0,T]})$.
This automatically implies that Assumption \ref{assume:psi} holds
for system (\ref{eq:system_spring}) with $\psi(\bfx,t)=x_1+\lambda
x_2$, $\lambda>0$. In particular, the following estimate holds
\[
\|\bfx(t)\|_{\infty,[t_0,T]}\leq
(1+\lambda^{-1})|x_1(t_0)|+(1+2\lambda^{-1})\|\psi(\bfx(t),t)\|_{\infty,[t_0,T]}
\]

Let us specify a class of control inputs $u$ which, in principle,
can ensure boundedness of solutions
$\bfx(t,\bfx_0,t_0,\thetavec,u)$ for every $\thetavec\in
\Omega_\theta$ and $\bfx_0\in\Real^n$. According to
(\ref{eq:assume_psi}), boundedness of
$\bfx(t,\bfx_0,t_0,\thetavec,u)$ is ensured if we find a control
input $u$ such that $\psi(\bfx(t),t)\in L_\infty^1[t_0,\infty]$.
To this objective consider the dynamics of system (\ref{system})
with respect to $\psi(\bfx,t)$:
\begin{eqnarray}\label{dpsi}
\dot{\psi}=L_{\bff(\bfx,\thetavec)}\psi(\bfx,t)+L_{\bfg(\bfx)}\psi(\bfx,t)u+\frac{\pd
\psi(\bfx,t)}{\pd t},
\end{eqnarray}
Assuming that the inverse
$\left(L_{\bfg(\bfx)}\psi(\bfx,t)\right)^{-1}$ exists everywhere,
we may choose the control input $u$ in the following class of
functions:
\begin{equation}\label{control}
\begin{split}
u(\bfx,\hat\thetavec,\omegavec,t)&=(L_{\bfg(\bfx)}\psi(\bfx,t))^{-1}\left(-L_{\bff(\bfx,\hat{\thetavec})}\psi(\bfx,t)-\varphi(\psi,\omegavec,t)-\frac{\pd\psi(\bfx,t)}{\pd
t}\right)\\
\varphi:& \ \Real\times\Real^w\times\Real_+\rightarrow\Real
\end{split}
\end{equation}
where $\omegavec\in\Omega_\omega\subset\Real^w$ is a vector of
{\it known} parameters of function $\varphi(\psi,\omegavec,t)$.
%Particular
%properties of the function $\varphi(\psi,\omegavec,t)$ in
%(\ref{control}) will be discussed later.
Denoting
$L_{\bff(\bfx,\thetavec)}\psi(\bfx,t)=f(\bfx,\thetavec,t)$ and
taking into account (\ref{control}) we may rewrite equation
(\ref{dpsi}) in the following manner:
\begin{eqnarray}\label{error_model}
\dpsi=f(\bfx,\thetavec,t)-f(\bfx,\hat{\thetavec},t)-\varphi(\psi,\omegavec,t)
\end{eqnarray}
Hence, feedback (\ref{control}) renders the original system
(\ref{system1}) into the well-known nonlinear {\it error model}
form \cite{Narendra89}\footnote{The error models
(\ref{error_model}) have proven to be convenient representations
of systems with nonlinear parametrization in the problems of
adaptation \cite{Annaswamy99,tpt2003_tac} and parameter estimation
\cite{Cao_2003}}.

In practical applications, state $\bfx$ of original model
(\ref{system1}) is hardly ever available. Furthermore, imprecise
physical models of the processes and measurement noise often lead
to the presence of {\it unmodeled dynamics} in
(\ref{error_model}). Although we do not address these issues in
detail in the present article, we do allow additive perturbations
that are functions of time from $L_{2}^1 [t_0,\infty]$ in the
right-hand side of (\ref{error_model}). In particular, instead of
(\ref{error_model}) we consider the following equation:
\begin{eqnarray}\label{error_model_d}
\dpsi=f(\bfx,\thetavec,t)-f(\bfx,\hat{\thetavec},t)-\varphi(\psi,\omegavec,t)+\varepsilon(t),
\end{eqnarray}
where, if not stated overwise, the function
$\varepsilon:\Real_+\rightarrow\Real$, $\varepsilon\in L_{2}^1
[t_0,\infty]\cap C^0$. One of the immediate advantages of
(\ref{error_model_d}) in comparison with (\ref{error_model}) is
that it allows us to take the presence of state observers in the
system into consideration. This  clearly widens the range of
possible applications of our results.% as discussed in Section 5.

% {\bf Can we provide any reference to
%your paper(s) discussing this?  Reviewers might ask for more
%details on this.}.

%which  due to the state observers or specific measurement noise.

%Function $\upsilon(t)$ in (\ref{error_model}), depending on the
%context, can serve as both additional control input or disturbance
%due to the model mismatches with the system dynamics.

Let us now specify the desired properties of function
$\varphi(\psi,\omegavec,t)$ in (\ref{control}),
(\ref{error_model_d}). The majority of known algorithms for
parameter estimation and adaptive control
\cite{Sastry89,Narendra89,Kokotovich95,Fradkov99} assume global
(Lyapunov) stability
 of system
(\ref{error_model_d}) for $\thetavec\equiv\hat{\thetavec}$.  In
our study, however, we would like to refrain from this standard
and at the same time restrictive requirement. Instead we propose
that finite energy of the signal
$f(\bfx(t),\thetavec,t)-f(\bfx(t),\hat{\thetavec}(t),t)$, defined
for example by its $L_{2}^1[t_0,\infty]$ norm with respect to the
variable $t$, results in finite deviation from the target set
given by equality $\psi(\bfx,t)=0$. Formally this requirement is
introduced in Assumption \ref{assume:gain}:
%{\bf This may be too restrictive. In this paper we propose that,
%instead of the Lyapunov stability, the finite energy of the signal
%$f(\bfx(t),\thetavec,t)-f(\bfx(t),\hat{\thetavec}(t),t)$, defined
%by its $L_{2}^1 [t_0,\infty]$ norm results in $\psi(\bfx,t)\in
%L_{\infty}^1$, as in Assumption \ref{assume:gain} below:}
\begin{assume}[Target dynamics operator]\label{assume:gain} Consider the following system:
\begin{equation}\label{eq:target_dynamics}
\dpsi=-\varphi(\psi,\omegavec,t)+\zeta(t),
\end{equation}
where $\zeta:\Real_+\rightarrow\Real$ and
$\varphi(\psi,\omegavec,t)$ is from (\ref{error_model_d}). Then
for every $\omegavec\in\Omega_\omega$ system
(\ref{eq:target_dynamics}) has $L_{2}^1 [t_0,\infty]\mapsto
L_\infty^1[t_0,\infty]$ gain with respect to input $\zeta(t)$. In
the other words,
\[
\zeta(t)\in L_{2}^1 [t_0,\infty]\Rightarrow
\psi(t,\psi_0,t_0,\omegavec)\in L_\infty^1[t_0,\infty], \
\psi_0\in\Real
\]
and there exists a function $\gamma_{\infty,2}$ such that
\begin{equation}\label{eq:gain_psi_L2}
\|\psi(t)\|_{\infty,[t_0,T]}\leq
\gamma_{\infty,2}(\psi_0,\omegavec,\|\zeta(t)\|_{2,[t_0,T]}), \ \
\forall \ \zeta(t)\in L_{2}^1[t_0,T]
\end{equation}
\end{assume}
In contrast to conventional approaches, Assumption
\ref{assume:gain} does not require global {\it asymptotic
stability} of the origin of (unperturbed, i.e for $\zeta(t)=0$)
system (\ref{eq:target_dynamics}). In fact, system
(\ref{eq:target_dynamics}) is allowed to have Lyapunov-unstable
equilibria. Moreover, there may be no equilibria  at all in
(\ref{eq:target_dynamics}), or it can even exhibit chaotic
dynamics. Examples of such systems, which potentially inherit
chaotic behavior but still satisfy Assumption \ref{assume:gain},
are the well-known Lorenz  \cite{Lorenz:1963} and Hindmarsh-Rose
\cite{Hindmarsh_and_Rose} oscillators.
%:
%\begin{equation}\label{eq:HR_model}
%\begin{split}
%    \dot{x}_{1} &= -a{x}_{1}^3+b{x}_{1}^2+x_{2}-x_{3}+I_0 + \zeta(t) \\
%    \dot{x}_{2} &= c-d{x}_{1}^2-x_{2}\\
%    \dot{x}_{3} &= \eps(s(x_{1}+x_0)-x_{3}),
%\end{split}
%\end{equation}
%where $a = 1$, $b = 3$, $c = 1$, $d = 5$, $s = 4$, $x_0 = 1.6$,
%$\eps = 0.001$, $I_0=1.4$.  System (\ref{eq:HR_model})
The last system models ion current through a membrane in the
living cell, and is widely used in artificial neural networks, for
instance, for processing of the visual information
\cite{Raffone_2003}.

When the stability of the target dynamics
$\dpsi=-\varphi(\psi,\omegavec,t)$ is known a-priori, one of the
benefits of Assumption \ref{assume:gain} is that there is no need
to know the {\it particular Lyapunov function} of the unperturbed
system. Apart from being, in some sense, a more friendly and less
invasive concept, this enables us to design adaptive/parameter
estimation procedures for systems with externally-driven
uncontrolled multistability
\cite{Arecchi:82,Weigel_98,Chizh_2000,Maldonado_1990}\footnote{In
systems with externally driven multistability, i.e. when there are
multiple coexistent attractors and trajectories switch from one
attractor to another depending on the external perturbation,
parameter estimation/control algorithms based on the knowledge of
a specific Lyapunov function require additional information about
the instant dynamical state (attractors and their allocation) of
the system itself. This leads to a necessity to {\it identify}
current dynamical state of the system {\it prior} to
control/identification of its parameters.}.

The differences between conventional restrictions on the goal
functionals and alternative requirements formulated in Assumptions
\ref{assume:psi}, \ref{assume:gain} are further illustrated with
Figure 1.
\begin{figure}\label{fig:shpere_bounds}
\begin{center}
\includegraphics[width=450pt]{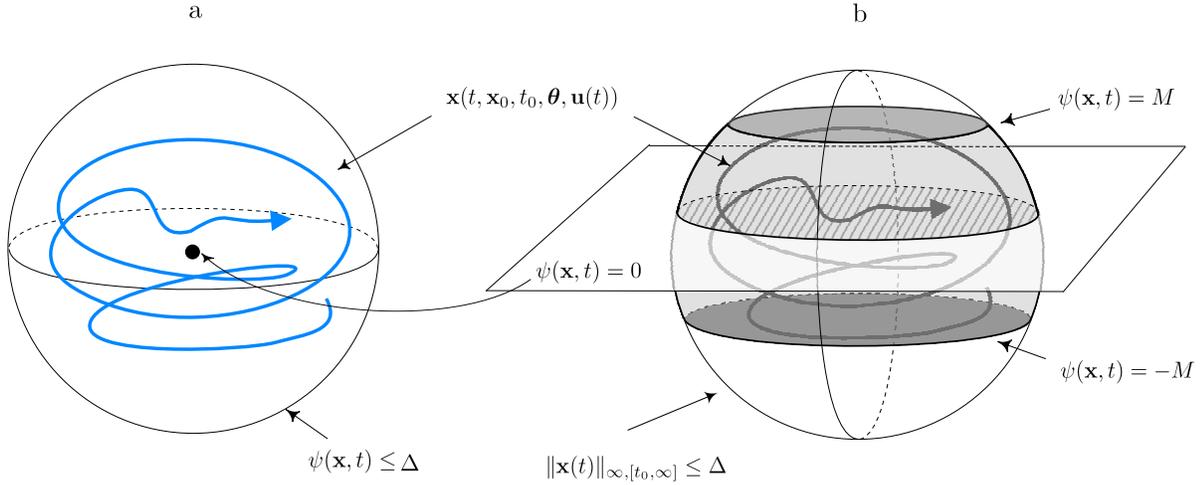}
\end{center}
\caption{Illustration to the choice of the goal functionals
$\psi(\bfx,t)$. Conventional requirements are illustrated with
panel a. Requirements following from Assumptions \ref{assume:psi},
\ref{assume:gain} are illustrated with panel b.}
\end{figure}
For simplicity it is assumed that function $\psi(\bfx,t)$ does not
depend on $t$ explicitly and therefore its zeroes form a (set)
surface in $\Real^n$. For the conventional approaches this set
should additionally satisfy metric conditions
(\ref{eq:standard_goal}) in $\Real^n$, Fig. 1. a. These conditions
often restrict class of the possible target sets to the points in
$\Real^n$. In case Assumptions \ref{assume:psi}, \ref{assume:gain}
are satisfied this restriction does not apply any more. Indeed,
given that $\zeta(t)\in L_2^1[t_0,\infty]$ we can bound
$\|\psi(\bfx(t),t)\|_{\infty,[t_0,\infty]}\leq
\gamma_{\infty,2}(\psi_0,\omegavec,\|\zeta(t)\|_{2,[t_0,\infty]})=M$.
Therefore, according to Assumption \ref{assume:gain}, the state
$\bfx(t)$ is bounded and belongs to the sphere
$\Omega_\bfx=\{\bfx\in\Real^n| \ \bfx: \ \|\bfx(t)\|\leq
\tilde{\gamma}(\bfx_0,\thetavec,M)=\Delta\}$. On the other hand,
the state $\bfx(t)$  belongs to the domain
$\Omega_{\psi}=\{\bfx\in\Real^n|\ \ \bfx: \ |\psi(\bfx,t)|\leq
M\}$. This implies that the segments of trajectory
$\bfx(t,\bfx_0,t_0,\thetavec,\bfu(t))$,  for $t\geq t_0$ will
remain in the  bounded domain $\Omega_\bfx\cap\Omega_\psi$
(shadowy volume in Fig. 1. b.) for all $t>t_0$.

The Figure 1 also emphasizes the difference between the proposed
operator framework and known approaches in adaptive control based
on geometrical representations \cite{Ortega2003}. Indeed, the
results based on coordinate transformations around the target
manifold (\ref{eq:target_set}) are applicable only in a subset of
$\Real^n$ where $\psi(\bfx,t)$ does not depend explicitly on $t$,
and rank of $\psi(\bfx,t)$ is constant. In this respect these
results are local. On the other hand, Assumptions
\ref{assume:psi}, \ref{assume:gain} do not require constant rank
conditions and allow both time-varying $\psi(\bfx,t)$ and
$\varphi(\psi,\omegavec,t)$. This makes Assumptions
\ref{assume:psi}, \ref{assume:gain} a suitable replacement to
conventional approaches for systems with non-stationary dynamics,
or ones which are far away from equilibrium or invariant target
manifolds.

%Another example is dynamics of multistable lasers
%\cite{Arecchi:82}
%\begin{equation}
%\begin{split}
%\dot{x} & = \tau^{-1}(y-k(t))x\\
%\dot{y} & = (y_0-y)\gamma - xy,\\
%k(t) & = k_0(1 + m \cos(2\pi f_d t)+\varepsilon(t)), \ k_0,m>0
%\end{split}
%\end{equation}
%where $\tau,\gamma>0$ are parameters, $k(t)$ is the total cavity
%losses, $f_d$ is the driving frequency.

So far we have introduced basic assumptions on system
(\ref{system1}) dynamics and the class of feedback considered in
this article. Let us now specify the class of functions
$f(\bfx,\thetavec,t)$ in (\ref{error_model_d}). Since general
parametrization of function $f(\bfx,\thetavec,t)$ is
methodologically difficult to deal with but solutions provided for
a restricted class of nonlinearities (for instance to those which
allow linear re-parametrization) often yield physically
implausible models, we have opted for a new class of
parameterizations. This class shall include a sufficiently broad
range of physical models, in particular those with nonlinear
parametrization; the proposed parameterizations will also, in
principle, be able to handle arbitrary state nonlinearity in the
class of functions from $\mathcal{C}^1$. As a candidate for such a
parametrization we suggest nonlinear functions that satisfy the
following assumption:
\begin{assume}[Monotonicity and Growth Rate in
Parameters]\label{assume:alpha}For the given function
$f(\bfx,\thetavec,t)$ in (\ref{error_model_d}) there exists
function  $\alphavec(\bfx,t): \Real^{n}\times \Real_+\rightarrow
\Real^d, \ \alphavec(\bfx,t)\in \mathcal{C}^1$ and positive
constant $D>0$ such that
\begin{equation}\label{eq:assume_alpha}
(f(\bfx,\hat{\thetavec},t)-f(\bfx,\thetavec,t))(\alphavec(\bfx,t)^{T}(\hat{\thetavec}-\thetavec))\geq0
\end{equation}
\begin{equation}\label{eq:assume_gamma}
|f(\bfx,\hat{\thetavec},t)-f(\bfx,\thetavec,t)|\leq D
|\alphavec(\bfx,t)^{T}(\hat{\thetavec}-\thetavec)|
\end{equation}
\end{assume}
The first inequality (\ref{eq:assume_alpha}) in Assumption
\ref{assume:alpha} holds, for example, for every smooth nonlinear
function which is monotonic with respect to a linear functional
$\phivec(\bfx)^{T}\thetavec$ over a vector of parameters:
\[
f(\bfx,\thetavec,t)=f_m(\bfx,\phivec(\bfx)^T\thetavec,t)
\]
\[
\sign\left( \frac{\pd f_m(\bfx,\lambda,t)}{\pd \lambda}\right) =\const %{\bf \{-1,0,1\} }
\]
Hence function $\alphavec(\bfx,t)$ satisfying
(\ref{eq:assume_alpha}) could be chosen in the following form:
$\alphavec(\bfx,t)=M \phivec(\bfx) \kappa(\bfx,t)$, where
$\kappa:\Real^n\times \Real_+\rightarrow \Real_+$,
$\kappa(\bfx,t)\in \mathcal{C}^1$.

The second inequality,
(\ref{eq:assume_gamma}), is satisfied if the function
$f(\bfx,\phivec(\bfx)^T\thetavec,t)$ does not grow faster than a
linear function in variable $\phivec(\bfx)^T\thetavec$ for every
$\bfx\in \Real^n$. This requirement holds, for example, for those
functions $f(\bfx,\phivec(\bfx)^T\thetavec,t)$ which are globally
Lipschitz in $\phivec(\bfx)^T\thetavec$:
\[
|f_m(\bfx,\phivec(\bfx)^T\thetavec,t)-f_m(\bfx,\phivec(\bfx)^T\thetavec',t)|\leq
D_{\theta}(\bfx,t)|\phivec(\bfx)^T(\thetavec-\thetavec')|
\]
In particular, inequalities (\ref{eq:assume_alpha}),
(\ref{eq:assume_gamma}) hold for the function
$f(\bfx,\phivec(\bfx)^T\thetavec,t)$ with $\alphavec(\bfx,t)=M
D_\theta(\bfx,t)\phivec(\bfx)$. A graphical illustration of the
choice of function $\alphavec(\bfx,t)$ is given in Figure 2.

\begin{figure}\label{fig:sector}
\begin{center}
\includegraphics[width=250pt]{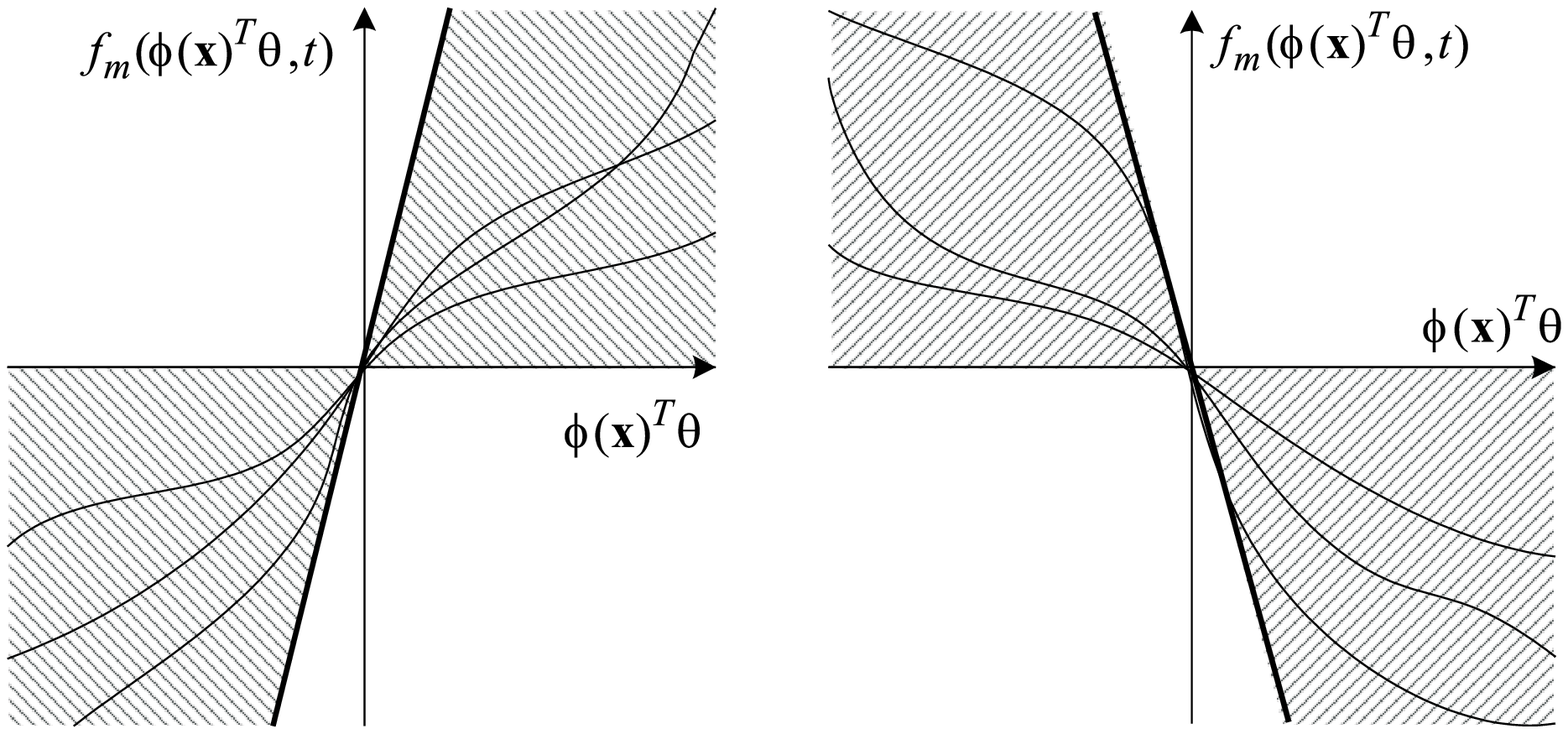}
\end{center}
\caption{Illustration to Assumption \ref{assume:alpha} for the
nonlinear in parameters function
$f(\bfx,\thetavec,t)=f_m(\phivec(\bfx)^{T}\thetavec,t)$,
$f_m:\Real\times \Real_+ \rightarrow \Real$. Thick lines stand for
the function $M D \phivec(\bfx)^T\thetavec$,
$D=\max_{\bfx,t}{|D_\theta(\bfx,t)|}$ in each block respectively.}
\end{figure}

This set of conditions naturally extends  from systems that are
linear in parameters  to those with nonlinear parametrization.
Assumption \ref{assume:alpha} covers (at least for bounded
$\thetavec,\hat{\thetavec}\in \Omega_\theta$) a considerable
variety of practically relevant models with nonlinear
parametrization. These include effects of stiction forces
\cite{Armstrong_1993}, slip and surface dependent friction given
by the ``magic formula" \cite{Pacejka91} or physics-inspired model
of the tyre \cite{Canudas_1999}, nonlinear processes in dampers
for automotive suspension \cite{Kitching_2000}, smooth saturation,
and dead-zones in mechanical systems. It further includes
nonlinearities in models of bio-reactors \cite{Boskovic_1995}. The
class of functions $f(\bfx,\thetavec,t)$ specified in Assumption
\ref{assume:alpha} can also serve as nonlinear replacement of
functions that are linear in their parameters in a variety of
piecewise approximation models. Last but not least, this set of
functions includes sigmoid and Gaussian nonlinearities, which are
favored in neuro and fuzzy control and mathematical models of
neural processes \cite{Abbott_2001}. Table 1 provides some of the
parametric nonlinearities that occur in these processes and their
corresponding functions $\alphavec(\bfx,t)$.
\begin{table}[htbp]\label{table:nonlinearities}
  \centering
  \caption{Examples of nonlinearities satisfying Assumption
  \ref{assume:alpha}. Parameters $\Delta_\theta$, $\Delta_r$ are positive constants}
\vskip 5mm {\small
  \begin{tabular*}{\textwidth}{@{\extracolsep{\fill}}|c|c|c|c|}
    \hline
     & & & \\
     physical meaning &  mathematical model & domain of &
    $\alphavec(\bfx,t)$\\
     & of uncertainty $f(\bfx,\thetavec,t)$ & physical&\\
     & &  relevance & \\
    \hline
     & & & \\
     stiction forces & $\theta_0 e^{-x_2^2 \theta_1}=e^{-x_2^2
\theta_1+\ln(\theta_0)}$ & $\Delta_\theta >\theta_0,\theta_1>0$ &$\alphavec(\bfx,t)=(-x_2^2,1)^T$ \\
   & $\bfx=(x_1,x_2)$ & $\bfx\in\Real^2$ & \\
  \hline
   & & & \\
     tyre-road friction \cite{Canudas_1999} & $F_n \sign(x_2)\frac{\frac{\sigma_0}{L}
G \frac{x_3}{1-x_3}}{\frac{\sigma_0}{L}
\frac{x_3}{1-x_3}+G} $& $\Delta_\theta>\theta>0$ & $\alphavec(\bfx,t)=\frac{x_3}{1-x_3}$ \\
 & $G =\theta\left(\mu_C +
(\mu_{S}-\mu_{C})e^{-\frac{|r x_2
x_3|}{|1-x_3|v_s}}\right)$ & $x_1,x_2\geq0$  & \\
  & $\bfx=(x_1,x_2,x_3)$ & $x_3\in(0,1)$ & \\
 & $F_n,\sigma_0,\mu_S,\mu_C>0$ - parameters &  & \\
 &  $\mu_S>\mu_C$ & & \\
\hline
 & & & \\
 force supported by & $\frac{K_o(\theta+1) A_p(x_1-x_2)}{L\left(\theta+K_o\frac{P_0}{x_3^2}\right)}$ & $\Delta_\theta>\theta>0$ & $\alphavec(\bfx,t)=A_p(x_1-x_2)$ \\
  hydraulic emulsion  & & $x_3>0$ &\\
  in suspension &$\bfx=(x_1,x_2,x_3)$ & & \\
   dampers \cite{Kitching_2000}  & $K_o,A_p,P_o,L>0$ - parameters& & \\
  & & & \\
 \hline
  & & & \\
 nonlinearities in & $\frac{x_1 x_2}{\theta_0 + \theta_1 x_1}, \ \frac{x_1 x_2}{\theta_0 + \theta_1 x_2}$ & $\Delta_\theta>\theta_0,\theta_1>0$ & $\alphavec(\bfx,t)=x_1 x_2(1, x_1)^T$ \\
 Monod's growth model & $\bfx=(x_1,x_2)$ & $x_1,x_2>0$ & $\alphavec(\bfx,t)=x_1 x_2(1, x_2)^T$\\
  of microorganisms \cite{Boskovic_1995}  & & & \\
 \hline
  & & & \\
  blur distortion model & $\sum_{i,j=1}^n e^{- \frac{|i - i_c|+|j-j_c|}{\theta}}r_{i,j}(t)$ & $\Delta_\theta>\theta>0$ & $\alphavec(\bfx,t)= 1$ \\
  in networks  & $r_{i,j}:\Real_+\rightarrow\Real_+$, $i,j,n\in \Natural$ & $\|r_{i,j}(t)\|_{\infty,[t_0,\infty]}\leq \Delta_r$ & \\
  for processing & $i_c,j_c\in \Natural$, $1\leq i_c,j_c\leq n$ & & \\
  of visual information & & & \\
 \hline
  \end{tabular*}
}
\end{table}

Assumption \ref{assume:alpha} bounds the growth rate of the
difference $|f(\bfx,\thetavec,t)-f(\bfx,\hat{\thetavec},t)|$ by
the functional
$D|\alphavec(\bfx,t)^{T}(\hat{\thetavec}-\thetavec)|$.  This will
help us to find a parameter estimation algorithm such that the
estimates converge to $\thetavec$ sufficiently fast for the
solutions of (\ref{system1}), (\ref{error_model_d})  to remain
bounded with non-dominating feedback (\ref{control}). On the other
hand, parametric error $\hat{\thetavec}-\thetavec$ can be inferred
from the changes in the variable $\psi(\bfx,t)$, according to
(\ref{error_model_d}), only by means of the difference
$f(\bfx,\thetavec,t)-f(\bfx,\hat{\thetavec},t)$. Therefore, as
long as convergence of the estimates $\hat{\thetavec}$ to
$\thetavec$ is expected, it is also useful to have  the estimate
of $|f(\bfx,\thetavec,t)-f(\bfx,\hat{\thetavec},t)|$ from below,
as specified in Assumption
\ref{assume:alpha_upper}\footnote{Despite Assumption
\ref{assume:alpha_upper} requires that
(\ref{eq:assume_alpha_upper}) holds for every $\bfx\in\Real^n$,
$\thetavec\in \Real^d$, and $t\in\Real_+$, we will see later that
for a variety of problems it is sufficient that it is satisfied
only locally.}:
\begin{assume}\label{assume:alpha_upper} For the given function
$f(\bfx,\thetavec,t)$ in (\ref{error_model_d}) and function
$\alphavec(\bfx,t)$, satisfying Assumption \ref{assume:alpha},
there exists a positive constant $D_1>0$  such that
\begin{equation}\label{eq:assume_alpha_upper}
|f(\bfx,\hat{\thetavec},t)-f(\bfx,\thetavec,t)|\geq D_1
|\alphavec(\bfx,t)^{T}(\hat{\thetavec}-\thetavec)|
\end{equation}
\end{assume}
% \footnote{As illustrated in  Table 1, both
%Assumptions \ref{assume:alpha}, \ref{assume:alpha_upper} hold at
%least locally for a variety of practically relevant
%nonlinearities}.

%In addition to the bounds provided in Assumptions
%\ref{assume:alpha}, \ref{assume:alpha_upper}, when analyzing
%performance of the closed loop system, it is also useful to have

In problems of parameter estimation, effectiveness of the
algorithms often depends on how "good" the nonlinearity
$f(\bfx,\thetavec,t)$ is, and how predictable locally is the
system's behavior. As a measure of goodness and predictability
usually the substitutes as smoothness, boundedness, and Lipschitz
conditions are considered. In our study, we distinguish several
such specific properties of functions $f(\bfx,\thetavec,t)$ and
$\varphi(\psi,\omegavec,t)$. These properties are provided below

\begin{hyp}\label{hyp:locally_bound_uniform_f} Function $f(\bfx,\thetavec,t)$ is locally bounded
with respect to $\bfx$, ${\thetavec}$  uniformly in $t$.
\end{hyp}

\begin{hyp}\label{hyp:locally_bound_uniform_df} Function $f(\bfx,\thetavec,t)\in \mathcal{C}^1$, and $ \pd
{f(\bfx,\thetavec,t)}/{\pd t}$ is locally bounded with respect to
$\bfx$, ${\thetavec}$ uniformly in $t$.
\end{hyp}

\begin{hyp}\label{hyp:globally_Lipschitz_uniform}
Function $f(\bfx,\thetavec,t)$ is globally Lipschitz with respect
to  $\thetavec$ uniformly in $\bfx$, $t$:
\[
\exists \ D_\theta>0:  \
|f(\bfx,\thetavec,t)-f(\bfx,\hat{\thetavec},t)|\leq
D_\theta\|\thetavec-\hat{\thetavec}\|
\]
\end{hyp}

\begin{hyp}\label{hyp:local_linear_bound} Let  $U_x\subset \Real^n$, $U_\theta\subset\Real^d$ be given and $U_x$, $U_\theta$ are
bounded. Then there exists constant $D_{U_x,U_\theta}>0$ such that
for every $\bfx\in U_x$ and $\thetavec,\hat{\thetavec}\in
U_\theta$ Assumption \ref{assume:alpha_upper} is satisfied with
$D_1=D_{U_x,U_\theta}$.
\end{hyp}

\begin{hyp}\label{hyp:locally_bound_uniform_phi} Function $\varphi(\psi,\omegavec,t)$ is locally bounded in $\psi$,
$\omegavec$ uniformly in $t$.
\end{hyp}

In the next section we present novel algorithms for adaptive
control and parameter estimation in nonlinear dynamical systems
(\ref{system}) which satisfy  Assumptions \ref{assume:psi},
\ref{assume:gain}, \ref{assume:alpha}, and
 \ref{assume:alpha_upper}. We show that under an additional
structural requirement, which relates properties of function
$\alphavec(\bfx,t)$ and vector-filed
$\bff(\bfx,\thetavec)=\bff_1(\bfx,\thetavec)\oplus\bff_2(\bfx,\thetavec)$
in (\ref{system1}), (\ref{system}), the following desired property
holds:
\begin{equation}\label{eq:desired_prop}
\bfx(t)\in L_\infty^n[t_0,\infty]; \
f(\bfx(t),\thetavec,t)-f(\bfx,\hat{\thetavec}(t),t)\in
L_{2}^1[t_0,\infty]
\end{equation}
After boundedness of the solutions is guaranteed, we prove that
\[
\lim_{t\rightarrow\infty}\psi(\bfx(t),t)=0
\]
In addition, we show that
\begin{equation}\label{eq:convergence}
\lim_{t\rightarrow\infty}\hat{\thetavec}(t)=\thetavec
\end{equation}
In particular we demonstrate that the standard persistent
excitation condition is sufficient to guarantee the convergence.
Furthermore, in the case that Assumptions
\ref{assume:alpha},\ref{assume:alpha_upper} hold only locally (for
$\bfx$ from a domain of $\Real^n$) we demonstrate that
sufficiently high excitation in the system still leads to the
desired estimates.

\section{Main Results}

Standard approaches  in parameter estimation and adaptation
problems usually assume feedback and a parameter adjustment
algorithm in the following form
\begin{equation}\label{eq:standard_control}
\begin{split}
& u=u(\bfx,\hat{\thetavec},t)\\
& \dot{\hat{\thetavec}}=\mathcal{A}_{\mathrm{lg}}(\psi,\bfx,t)
\end{split}
\end{equation}
The most favorite strategy of finding these, known also as {\it
certainty-equivalence principle}, is a two-stage design
prescription.  First, construct uncertainty-dependent feedback
$u(\bfx,\thetavec,t)$, $\thetavec\in\Omega_\theta$ which ensures
boundedness of the trajectories $\bfx(t)$. Second, replace
$\thetavec$ with $\hat{\thetavec}$ in $u(\bfx,\thetavec,t)$ and,
given the constraints (e.g., $\dot{\bfx}$, $\thetavec$ cannot be
measured explicitly, while state $\bfx$ is available), design
function
$\mathcal{A}_{\mathrm{lg}}:\Real\times\Real^n\times\Real_+\rightarrow\Real^d$
which guarantees (\ref{eq:desired_prop}), (\ref{eq:convergence}),
or/and $\psi(\bfx(t),t)\rightarrow 0$.

With this strategy, the design of the feedback
$u(\bfx,\thetavec,t)$
%can be done independently to a certain
%extent
is generally independent\footnote{In particular, it is the
standard requirement that function  $u(\bfx,\thetavec,t)$ should
guarantee Lyaponov stability of the system for
$\hat{\thetavec}=\thetavec$, while parameter adjustment algorithms
use this property in order to ensure stability of the whole
system. No other properties are required from the function
$u(\bfx,\thetavec,t)$.} of the specific design of the parameter
estimation algorithm $\mathcal{A}_\mathrm{lg}(\psi,\bfx,t)$.
% on the parameter estimation algorithm
%$\mathcal{A}_\mathrm{lg}(\psi,\bfx,t)$.
This allows the full benefit of
%{\bf Such decoupling utilizes the benefits of
contemporary nonlinear control theory
\cite{Isidory_96,Isidory_99,Khalil_2002,Nijmeijer_90}  in
designing feedbacks $\bfu(\bfx,\thetavec,t)$. On the  other hand,
this strategy equally benefits from conventional parameter
estimation and adaptation theories
\cite{Ljung_99,Eykhoff,Narendra89,Fradkov99} which provide a list
of the ready-to-be-implemented algorithms under the {\it
assumption} that feedback $\bfu(\bfx,\thetavec,t)$ ensures
stability of the system.

%the conventional parameter estimation and adaptation theory
%\cite{Ljung_99,Eykhoff,Narendra89,Fradkov99}, when designing
%feedback $u(\bfx,\thetavec,t)$.}

% under assumptions that
%$\thetavec$ are known.

%{\bf In general, the certainty-equivalence principle does not take
%into account possible interactions between stabilizing control and
%parameter estimation procedures.

Ironically, the power of the certainty-equivalence principle --
simplicity and relative independence of the stages of design -- is
also its Achilles' heel. This principle does not take into account
the  possible interactions between stabilizing control and
parameter estimation procedures. It has been reported in
\cite{Stotsky93,Ortega02,Ortega2003,tpt2003_tac} that an
additional ``interaction" term
$\hat{\thetavec}_P(\bfx,t):\Real^n\times\Real_+\rightarrow
\Real^d$ added to the parameters ${\thetavec}$ in function
$u(\bfx,{\thetavec},t)$:
  $u(\bfx,\thetavec +
\hat{\thetavec}_P(\bfx,t),t)$ introduces new properties to the
system. Unfortunately, straightforward introduction of this
"interaction" term as a new variable of the design affects its
simplicity, internal order, and so much favored independence of
the design stages (control and estimation).

An alternative strategy which introduces a new design paradigm is
proposed in \cite{t_fin_formsA&T2003,ECC_2003}. Its main idea is
that adaptation algorithms in (\ref{eq:standard_control}) are
initially allowed to depend on unmeasurable variables
$\dot{\psi},\dot{\bfx},\thetavec$
\begin{equation}\label{eq:ideal_alg}
\dot{\hat{\thetavec}}=\mathcal{A}_{\mathrm{lg}}^\ast(\psi,\dpsi,\bfx,\dot{\bfx},\thetavec,t)
\end{equation}
For this reason we refer to such algorithms as {\it virtual
algorithms}. If the desired properties (\ref{eq:desired_prop}),
(\ref{eq:convergence}) are ensured with (\ref{eq:ideal_alg}) then,
taking into account properties of the vector-fields
$\bff(\bfx,\thetavec)$, $\bfg(\bfx)$ in (\ref{system1}), we
convert the unrealizable algorithm (\ref{eq:ideal_alg}) into an
equivalent representation in integro-differential, or finite, form
\cite{Fradkov86}:
\begin{equation}\label{eq:new_control}
\begin{split}
& \hat{\thetavec}=\Gamma(\hat{\thetavec}_P(\bfx,t)+\hat{\thetavec}_{I}(t)), \ \Gamma\in\Real^{d\times d},  \ \Gamma>0\\
& \dot{\hat{\thetavec}}_I=\mathcal{A}_{\mathrm{lg}}(\psi,\bfx,t)
\end{split}
\end{equation}
This  approach preserves the convenience of the
certainty-equivalence principle, as the feedback
$u(\bfx,\thetavec,t)$ could, in principle, be built independently
of the subsequent parameter adjustment procedure. At the same
time, it provides the necessary interaction term
$\hat{\thetavec}_P(\bfx,t)$ ensuring the required properties
(\ref{eq:desired_prop}), (\ref{eq:convergence}) of the closed-loop
system even if function $f(\bfx,\thetavec,t)$ in
(\ref{error_model_d}) is nonlinear in $\thetavec$.

%A diagram, illustrating
%this strategy is provided below:
%\begin{equation}\label{eq:diagram}
%\begin{array}{llrll}
%\mathrm{Goal}&\rightarrow &
%\dot{\hat{\thetavec}}=\mathcal{A}_{\mathrm{lg}}^\ast(\psi,\dpsi,\bfx,\dot{\bfx},\thetavec,t)&
%\rightarrow &  \hat{\thetavec}(\bfx,t)=\hat{\thetavec}_P(\bfx,t)+\int_0^{t}\mathcal{A}_\mathrm{lg}(\psi(\bfx(\tau),\tau),\bfx(\tau),\tau)d\tau\\
%& & &  \nearrow &\\
%& &\mathrm{Properties \ of \ } \bff(\bfx,\thetavec,t)& &
%\end{array}
%\end{equation}

In this paper we propose the following class of {\it virtual
adaptation algorithms}\footnote{Choice of the virtual algorithm in
the form of equation (\ref{eq:virtual_alg}) is motivated by our
previous study of derivative-dependent algorithms for systems with
uncertainties that are nonlinear in their parameters
\cite{tpt2002_at,tpt2003_tac}}:
\begin{equation}\label{eq:virtual_alg}
\dot{\hat{\thetavec}}=\Gamma(\dpsi+\varphi(\psi,\omegavec,t))\alphavec(\bfx,t)
+ \mathcal{Q}(\bfx,\hat{\thetavec},t)(\thetavec-\hat{\thetavec}),
\ \Gamma\in \Real^{d\times d}, \  \Gamma>0
\end{equation}
where
$\mathcal{Q}(\bfx,\hat{\thetavec},t):\Real^n\times\Real^d\times\Real_+\rightarrow
\Real^{d\times d}$, $\mathcal{Q}(\cdot)\in \mathcal{C}^0$. As a
candidate for finite form realization (\ref{eq:new_control}) of
algorithms (\ref{eq:virtual_alg}) we select the following set of
equations:
\begin{equation}\label{fin_forms_ours_tr1}
\begin{split}
\hat{\thetavec}(\bfx,t)&=\Gamma(\hat{\thetavec}_P(\bfx,t)+\hat{\thetavec}_I(t));
\  \Gamma\in\Real^{d\times d}, \ \Gamma>0
\\ \hat{\thetavec}_P(\bfx,t)&=
\psi(\bfx,t)\alphavec(\bfx,t)-\Psi(\bfx,t) \\
\dot{\hat{\thetavec}}_I&=\varphi(\psi(\bfx,t),\omegavec,t)\alphavec(\bfx,t)+\mathcal{R}(\bfx,\hat{\thetavec},u(\bfx,\hat{\thetavec},t),t),
\end{split}
\end{equation}
where function
$\Psi(\bfx,t):\Real^{n}\times\Real_+\rightarrow\Real_d$,
$\Psi(\bfx,t)\in \mathcal{C}^1$  satisfies Assumption
\ref{assume:explicit_realizability}
\begin{assume}\label{assume:explicit_realizability} There exists function $\Psi(\bfx,t)$ such that
\begin{equation}\label{eq:assume_explicit}
\frac{\pd \Psi(\bfx,t)}{\pd \bfx_2}-\psi(\bfx,t)\frac{\pd
\alphavec(\bfx,t)}{\pd \bfx_2}=\mathcal{B}(\bfx,t),
\end{equation}
where
$\mathcal{B}(\bfx,t):\Real^n\times\Real_+\rightarrow\Real^{d\times
p}$ is either zero or, if $\bff_2(\bfx,\thetavec)$ is
differentiable in $\thetavec$, satisfies the following:
\[
\mathcal{B}(\bfx,t)\mathcal{F(\bfx,\thetavec,\thetavec')}\geq 0 \
\ \forall \thetavec,\thetavec'\in \Omega_\theta, \ \bfx\in \Real^n
\]
\[
 \mathcal{F(\bfx,\thetavec,\thetavec')}=\int_0^1 \frac{\pd
\bff_2(\bfx, \bfs(\lambda))}{\pd \bfs} d\lambda, \ \
\bfs(\lambda)=\thetavec'\lambda+\thetavec(1-\lambda)
\]
\end{assume}
Function
$\mathcal{R}(\bfx,\hat{\thetavec},u(\bfx,\hat{\thetavec},t),t):\Real^n\times\Real^d\times\Real\times\Real_+\rightarrow\Real^d$
in (\ref{fin_forms_ours_tr1}) is given as follows:
\begin{eqnarray}\label{fin_forms_ours_tr11}
\mathcal{R}(\bfx,u(\bfx,\hat{\thetavec},t),t)&=&{\pd
\Psi(\bfx,t)}/{\pd t}-\psi(\bfx,t)({\pd
\alphavec(\bfx,t)}/{\pd t})-\nonumber\\
& &  (\psi(\bfx,t)L_{\bff_1} \alphavec(\bfx,t)-L_{\bff_1}
\Psi(\bfx,t))-(\psi(\bfx,t)L_{\bfg_1}\alphavec(\bfx,t)-L_{\bfg_1}
\Psi(\bfx,t))u(\bfx,\hat{\thetavec},t)\nonumber \\
& & +
\mathcal{B}(\bfx,t)(\bff_2(\bfx,\hat{\thetavec})+\bfg_2(\bfx)u(\bfx,\hat{\thetavec},t)).
\end{eqnarray}
Functions $\Psi(\bfx,t)$ and
$\mathcal{R}(\bfx,\hat{\thetavec},u(\bfx,\hat{\thetavec},t),t)$
are introduced into (\ref{fin_forms_ours_tr1}) in order to shape
the derivative $\dot{\hat{\thetavec}}(\bfx,t)$ to fit equation
(\ref{eq:virtual_alg}).  The role of function $\Psi(\bfx,t)$ in
(\ref{fin_forms_ours_tr1}) is to compensate for the
uncertainty-dependent term
$\psi(\bfx,t)L_{\bff_2(\bfx,\thetavec)}\alphavec(\bfx,t)$, and
equation (\ref{eq:assume_explicit}) is the condition when such
compensation is possible\footnote{We show below, in the proof of
Theorem \ref{stability_theorem} (see Appendix), that Assumption
\ref{assume:explicit_realizability} is indeed sufficient for the
function $\hat{\thetavec}(\bfx,t)$ to be a realization of
(\ref{eq:virtual_alg}).}.  With the function
$\mathcal{R}(\bfx,\hat{\thetavec},u(\bfx,\hat{\thetavec},t),t)$ we
eliminate the  influence of the uncertainty-independent vector
fields $\bff_1(\bfx)$, $\bfg_1(\bfx)$, and $\bfg_2(\bfx)$ on the
desired form of the time-derivative
$\dot{\hat{\thetavec}}(\bfx,t)$.
%{\bf Equation (20) is still not well motivated.  Why not to explain every term in (20)?
%Further, you seem to use $L_{\bff(\bfx,\thetavec)}\psi(\bfx,t)=f(\bfx,\thetavec,t)$ in Theorem 1,
%but in (20) you use $L_{\bff_1}$ and $L_{\bff_2}$.  Possible confusion?}
The properties of system (\ref{system1}), together with control
(\ref{control}) and this new adaptation algorithm
(\ref{fin_forms_ours_tr1}), (\ref{fin_forms_ours_tr11}), are
summarized in Theorem \ref{stability_theorem} and Theorem
\ref{convergence_theorem}.

\begin{thm}[Boundedness]\label{stability_theorem}
Let system (\ref{system1}), (\ref{error_model_d}),
(\ref{fin_forms_ours_tr1}), (\ref{fin_forms_ours_tr11}) be given
and Assumptions
\ref{assume:alpha},\ref{assume:alpha_upper},\ref{assume:explicit_realizability}
be satisfied. Then the following properties hold

P1) Let for the given initial conditions $\bfx(t_0)$,
$\hat{\thetavec}_I(t_0)$ and parameters vector $\thetavec$,
interval $[t_0,T^\ast]$ be the (maximal) time-interval of
existence of solutions in the closed loop system (\ref{system1}),
(\ref{error_model_d}), (\ref{fin_forms_ours_tr1}),
(\ref{fin_forms_ours_tr11}). Then
\[
f(\bfx(t),\thetavec,t)-f(\bfx(t),\hat{\thetavec}(t),t))\in L_{2}^1
[t_0,T^\ast]
\]
and
\begin{equation}\label{eq:f_diff_L2}
\begin{split}
&
\|f(\bfx(t),\thetavec,t)-f(\bfx(t),\hat{\thetavec}(t),t))\|_{2,[t_0,T^\ast]}\leq
D_f(\thetavec,t_0,\Gamma,\|\varepsilon(t)\|_{2,[t_0,T^\ast]});\\
&
D_f(\thetavec,t_0,\Gamma,\|\varepsilon(t)\|_{2,[t_0,T^\ast]})=\left(\frac{D}{2}\|\thetavec-\hat{\thetavec}(t_0)\|^{2}_{\Gamma^{-1}}\right)^{0.5}
+ \frac{D}{D_1}\|\varepsilon(t)\|_{2,[t_0,T^\ast]}  \\
& \|\thetavec-\hat\thetavec(t)\|^{2}_{\Gamma^{-1}}\leq
\|\hat{\thetavec}(t_0)-\thetavec\|^{2}_{\Gamma^{-1}}+\frac{D}{2
D_1^2}\|\varepsilon(t)\|^{2}_{2,[t_0,T^\ast]}
\end{split}
\end{equation}

\noindent  In addition, if Assumptions \ref{assume:psi} and
\ref{assume:gain} are satisfied then

P2)  $\psi(\bfx(t),t)\in L_\infty^1[t_0,\infty]$, $\bfx(t)\in
L_{\infty}^n[t_0,\infty]$ and
\begin{equation}\label{eq:psi_gain}
\|\psi(\bfx(t),t)\|_{\infty,[t_0,\infty]}\leq
\gamma_{\infty,2}\left(\psi(\bfx_0,t_0),\omegavec,D_f(\thetavec,t_0,\Gamma,\|\varepsilon(t)\|_{2,[t_0,\infty]})+\|\varepsilon(t)\|_{2,[t_0,\infty]}\right)
\end{equation}

%{\bf One minor and one major concern.  Minor: $\|\bfx\|_{\Gamma}^2$ in Section 2, but here $\|\bfx\|_{\Gamma^{-1}}^2$.

%Major: Your P1 and P2 are likely to turn the reviewer off because they appear too complex.  Can we state that
%the right-hand sides are bounded, i.e., $D_f(\cdot)$ and $\gamma$?  It's not important here what particular
%form those bounds may be (whoever is interested would look into the proof), so stating that they are bounded
%is both sufficient and gives us a chance not to loose the reviewer.  Realistically, it's very unlikely that anybody
%will compute your bounds anyway because they contain hard-to-find terms.
%}

P3) if properties  H\ref{hyp:locally_bound_uniform_f},
H\ref{hyp:locally_bound_uniform_phi} hold, and system
(\ref{eq:target_dynamics}) has $L_{2}^1 [t_0,\infty]\mapsto
L_{p}^1 [t_0,\infty]$, $p>1$ gain with respect to input $\zeta(t)$
and output $\psi$ then
\begin{equation}\label{eq:convergence_psi_theorem}
\varepsilon(t)\in L_{2}^1 [t_0,\infty]\cap
L_{\infty}^1[t_0,\infty]\Rightarrow
\lim_{t\rightarrow\infty}\psi(\bfx(t),t)=0
\end{equation}

If, in addition, property H\ref{hyp:locally_bound_uniform_df}
holds, and functions $\alphavec(\bfx,t)$, $\pd \psi(\bfx,t)/\pd t$
are locally bounded with respect to $\bfx$ uniformly in $t$, then

P4) the following limiting relation holds
\begin{equation}\label{eq:convergence_f_theorem}
\lim_{t\rightarrow\infty}f(\bfx(t),\thetavec,t)-f(\bfx(t),\hat{\thetavec}(t),t)=0
\end{equation}

\end{thm}
Proofs of Theorem \ref{stability_theorem} and subsequent results
are given in the Appendix.

Before we proceed with discussion of the results of Theorem
\ref{stability_theorem}, we wish to comment on Assumption
\ref{assume:explicit_realizability}.
%This assumption provides
%specific restrictions on the choice of function $\Psi(\bfx,t)$
%which allow derivative-independent realization
%(\ref{eq:new_control}) of virtual algorithm
%(\ref{eq:virtual_alg}). As soon as this step is done, then
%properties P1)--P4) of the Theorem follow from the properties of
%the virtual algorithm. In particular,
Assumption
\ref{assume:explicit_realizability} links  the possibility to
design the parameter adjustment algorithm in the form of equation
(\ref{eq:virtual_alg}),
%which enables ability to deal with nonlinear parametrization \cite{tpt2002_at,tpt2003_tac},
%{\bf This statement is removed because virtual algorithm by itself is not an enabler of anything.
%Further, mentioning nonlinear parametrization is a distraction here.}
with the properties of functions $\alphavec(\bfx,t)$ and
$\psi(\bfx,t)$. These functions  depend on the properties of
nonlinearity $f(\bfx,\thetavec,t)$ itself (function
$\alphavec(\bfx,t)$) and, importantly, on the chosen specification
of the desired target set:
\[
\{\bfx\in\Real^n|\psi(\bfx,t)=0\}\subseteq\Omega_0
\]
given by function $\psi(\bfx,t)$. Specific properties of the
functions $f(\bfx,\thetavec,t)$ and $\psi(\bfx,t)$ are
interrelated through the possibility to solve partial differential
equation (\ref{eq:assume_explicit}) for the function
$\Psi(\bfx,t)$.
 Let
$\mathcal{B}(\bfx,t)=\mathrm{col}(\mathcal{B}_1(\bfx,t)$,$\dots$,
$\mathcal{B}_d(\bfx,t))$, and $\alphavec(\bfx,t)\in
\mathcal{C}^2$,
$\alphavec(\bfx,t)=\mathrm{col}(\alpha_1(\bfx,t),\dots,\alpha_d(\bfx,t))$,
then necessary  and sufficient conditions for existence of the
function $\Psi(\bfx,t)$ follow from the
Poincar$\acute{\mathrm{e}}$ lemma:
\begin{equation}\label{eq:poincare}
\frac{\pd}{\pd \bfx_2}\left(\psi(\bfx,t)\frac{\pd
\alpha_i(\bfx,t)}{\pd
\bfx_2}+\mathcal{B}_i(\bfx,t)\right)=\left(\frac{\pd}{\pd
\bfx_2}\left(\psi(\bfx,t)\frac{\pd \alpha_i(\bfx,t)}{\pd
\bfx_2}+\mathcal{B}_i(\bfx,t)\right)\right)^T
\end{equation}
This relation, in the form of conditions of existence of the
solutions for function $\Psi(\bfx,t)$ in
(\ref{eq:assume_explicit}), takes into account structural
properties of system (\ref{system1}), (\ref{error_model_d}).
Indeed, let $\mathcal{B}(\bfx,t)=0$ and consider partial
derivatives $\pd \alpha_i(\bfx,t)/\pd \bfx_2$, $\pd
\psi(\bfx,t)/\pd \bfx_2$ with respect to vector
$\bfx_2=(x_{21},\dots,x_{2p})^T$. Let
\begin{equation}\label{eq:single_dim}
\begin{split}
\frac{\pd \psi(\bfx,t)}{\pd \bfx_2}&=\left(\begin{array}{cccccccc}
                                     0& 0
                                     & \cdots & 0& \ast & 0&\cdots&0
                                    \end{array}\right)\\
\frac{\pd
\alpha_i(\bfx,t)}{\pd\bfx_2}&=\left(\begin{array}{cccccccc}
                                     0 & 0
                                     & \cdots & 0&
                                     \ast &
                                     0&\cdots&0
                                        \end{array}\right)
\end{split}
\end{equation}
where symbol $\ast$ denotes a function of $\bfx$ and $t$. Then
condition (\ref{eq:single_dim}) guarantees that equality
(\ref{eq:poincare}) (and, subsequently, Assumption
\ref{assume:explicit_realizability}) holds. Whether or not
Assumption \ref{assume:explicit_realizability}  holds, depends,
roughly speaking, on how large is the part of partition $\bfx_2$
that enters the arguments of functions $\psi(\bfx,t)$,
$\alphavec(\bfx,t)$.
%In other words, how high the dimensionally of
%the uncertainty $f(\bfx,\thetavec,t)$ and function $\psi(\bfx,t)$
%are connected with the uncertainty-dependent partition of system
%(\ref{system1}) itself.
%{\bf The removed sentence is not clear, and it'd add nothing anyway.}
In the case of $\pd \alpha(\bfx_1\oplus \bfx_2,t)/\pd \bfx_2=0$,
Assumption \ref{assume:explicit_realizability} holds for arbitrary
$\psi(\bfx,t)\in \mathcal{C}^1$. If $\psi(\bfx,t)$,
$\alphavec(\bfx,t)$ depend on just a single component of $\bfx_2$,
for instance $x_{2k}, \ k\in\{0,\dots,p\}$, then conditions
(\ref{eq:single_dim}) hold and function $\Psi(\bfx,t)$ can be
derived explicitly by integration
\begin{equation}\label{eq:single_dim_int}
\Psi(\bfx,t)=\int\psi(\bfx,t)\frac{\alphavec(\bfx,t)}{\pd x_{2k}}d
x_{2k}
\end{equation}
In all other cases, the existence of the required function
$\Psi(\bfx,t)$ follows from (\ref{eq:poincare}).

The necessity to satisfy Assumption
\ref{assume:explicit_realizability} may seem to be a critical
restriction, which limits  applicability of our approach. However,
we notice that it holds in the relevant problem
settings\footnote{See, for example, the problem setting  in
\cite{Cao_2003} for parameter estimation in the presence of
nonlinear state-dependent parametrization. This problem setting,
according to our knowledge, is by far one of the most general
available in the literature.} for arbitrary $\alphavec(\bfx,t),
\psi(\bfx,t)\in \mathcal{C}^1$. Consider, for instance
\cite{Cao_2003}, where the class of systems is restricted to
(\ref{eq:cao_annaswamy_setup}):
\begin{equation}\label{eq:cao_annaswamy_setup}
\dot{x}=-\varrho(x,u)x+f(\thetavec,u,x), \
\varrho(x,t)>\varrho_{\min}>0, \ x\in \Real
\end{equation}
The dimension of the state in system
(\ref{eq:cao_annaswamy_setup}) coincides with that of the
uncertainty-dependent partition and equals to unit
($\dim\{\bfx\}=\dim\{\bfx_2\}=1$). Hence, according to
(\ref{eq:single_dim_int}), and in case functions $\psi(x,t), \
\alphavec(x,t)\in \mathcal{C}^1$, there will always exist a
function $\Psi(x,t)$ satisfying equality
(\ref{eq:assume_explicit}) with $\mathcal{B}(x,t)=0$.

In the general case, when $\dim\{\bfx_2\}>1$, the problems of
finding function $\Psi(\bfx,t)$ satisfying condition
(\ref{eq:assume_explicit}) can be avoided (or converted into one
with an already known solutions such as (\ref{eq:poincare}),
(\ref{eq:single_dim_int})) by the {\it embedding} technique
proposed in \cite{ECC_2003}. The main idea of the method is to
introduce an auxiliary system
\begin{equation}\label{eq:embed}
\begin{split}
\dot{\xivec}&=\bff_\xivec(\bfx,\xivec,t), \ \xivec\in\Real^z \\
\bfh_\xi&=\bfh_\xi(\xivec,t), \
\Real^z\times\Real_+\rightarrow\Real^h
\end{split}
\end{equation}
such that
\begin{equation}\label{eq:embed_L2}
f(\bfx(t),\thetavec,t)-f(\bfx_1(t)\oplus\bfh_\xi(t)\oplus\bfx_2'(t),\thetavec,t)\in
L_{2}^1 [t_0,\infty]
\end{equation}
and $\dim\{{\bfh_\xi}\}+\dim{\{\bfx_2'\}}=p$. Then
(\ref{error_model_d}) can be rewritten as follows:
\begin{eqnarray}\label{error_model_d1}
\dpsi=f(\bfx_1\oplus\bfh_\xi\oplus\bfx_2',\thetavec,t)-f(\bfx_1\oplus\bfh_\xi\oplus\bfx_2',\hat\thetavec,t)-\varphi(\psi,\omegavec,t)+\varepsilon_\xi(t),
\end{eqnarray}
where $\varepsilon_\xi(t)\in L_{2}^1 [t_0,\infty]$, and
$\dim\{\bfx_2'\}=p-h<p$. In principle, the dimension of $\bfx_2'$
could be reduced to $1$ or $0$. As soon as this is ensured,
Assumption \ref{assume:explicit_realizability} will be satisfied
and the results of Theorem \ref{stability_theorem} follow.
Sufficient conditions ensuring the existence of such an embedding
in general case are provided in \cite{ECC_2003}. For systems in
which the parametric uncertainty  can be reduced to vector fields
with low-triangular structure the embedding is given in
\cite{ALCOSP_2004}.

An alternative way to construct system (\ref{eq:embed}) with the
desired properties is to use (possible, high-gain, discontinuous)
robust observers. In order to illustrate this approach consider
the rather general case when function $\bff_2(\bfx,\thetavec)$ in
(\ref{system1}) is given as
$\bff_2(\bfx,\thetavec)=\bar{\bff}_2(\bfx)+\phivec(\bfx,\thetavec)$,
and function $\phivec(\bfx,\thetavec)$ is bounded.

Let in addition there exist continuous functions
$\bfh_\epsilon:\Real^p\rightarrow\Real^p$,
$\bfh_\xi:\Real^p\rightarrow\Real^p$ such that the following
inequality is satisfied
\begin{equation}\label{eq:embed_bound_for_L2}
\|\bfh_\epsilon(\xivec-\bfx_2)\|\geq
|f(\bfx_1\oplus\bfh_\xi(\xivec),\thetavec,t)-f(\bfx_1\oplus\bfx_2,\thetavec,t)|
\end{equation}
As a candidate for yet unknown tracking system (\ref{eq:embed}) we
select the following
\begin{equation}\label{eq:embed_aux}
\dot{\xivec}=\bar{\bff}(\bfx)+\bff_{\epsilon}(\xivec-\bfx_2)+\bfg_2(\bfx)u
+ \boldsymbol{\upsilon}
\end{equation}
where function $\bff_{\epsilon}:\Real^p\rightarrow\Real^p$ and
auxiliary input $\boldsymbol{\upsilon}\in\Real^p$  are the design
parameters. Subtracting equations for $\dot{\bfx}_2$ in
(\ref{system1}) from (\ref{eq:embed_aux})  yields:
\begin{equation}\label{eq:embed_error_dyn}
\begin{split}
\dot{\epsvec}&=\bff_{\epsilon}(\epsvec)-\phivec(\bfx(t),\thetavec)+\boldsymbol{\upsilon}\\
\bfy_\epsilon&=\bfh_\epsilon(\epsvec)
\end{split}
\end{equation}
where $\epsvec=\xivec-\bfx_2$. Let us finally choose the function
$\bff_\epsilon$ in (\ref{eq:embed_error_dyn}) such that the system
$\dot{\epsvec}=\bff_{\epsilon}(\epsvec)+\boldsymbol{\upsilon}$ is
strictly passive with a positive definite storage function
$V(\epsvec,t):$
\begin{equation}\label{eq:storage}
\dot{V}(\epsvec,t)\leq \bfy_\epsilon^T \boldsymbol{\upsilon} -
\beta \|\bfy_\epsilon\|^2, \ \beta>0
\end{equation}
According to \cite{Lorea_2001}\footnote{In \cite{Lorea_2001} one
extra assumption on the function $\bff_\epsilon(\epsvec)$ in
(\ref{eq:embed_error_dyn}) is imposed. In particular it is
required that the system
$\dot{\epsvec}=\bff_{\epsilon}(\epsvec)+\upsilon$ is strongly
zero-detectable with respect to inputs $\upsilon$ and output
$\bfy_\epsilon$. In our case, however, the limiting relations
$\lim_{t\rightarrow\infty}\bfy_\epsilon(t) =0$,
$\lim_{t\rightarrow\infty}\epsvec(t) =0$ are not necessary.
Therefore, as follows from the proof of Theorem 2 in
\cite{Lorea_2001}, in order to show just
$\|\bfy_\varepsilon(t)\|\in L_2^p[t_0,\infty]$ the assumption of
strong zero-delectability can be omitted.}(page 1484, Theorem 2)
inequality  (\ref{eq:storage}) guarantees that there always exists
input $\boldsymbol{\upsilon}(t)$ in (\ref{eq:embed_aux}) such that
$\bfh_\epsilon(\epsvec(t))\in L_{2}^p [t_0,\infty]$. Then taking
into account (\ref{eq:embed_bound_for_L2}) we can conclude that
condition (\ref{eq:embed_L2}) holds (with $\bfx_2': \
\dim\{\bfx_2'\}=0$). This implies that the original error model
(\ref{error_model_d}) can be converted into
(\ref{error_model_d1}), which in our case satisfies Assumption
\ref{assume:explicit_realizability} (${\pd \alpha_i(\bfx,t)}{\pd
\bfx_2}=0$ for the corresponding $\alpha_i(\bfx,t)$ in
(\ref{error_model_d1})).

Let us now briefly comment on the results of Theorem
\ref{stability_theorem}. The theorem ensures a set of relevant
properties for both control (P2, P3) and parameter estimation
problems (P1, P4). These properties, as illustrated with
(\ref{eq:f_diff_L2})--(\ref{eq:convergence_f_theorem}), provide
conditions for boundedness of the solutions
$\bfx(t,\bfx_0,t_0,\thetavec,u(t))$, reaching the target set
$\Omega_0$, and exact compensation of the uncertainty term
$f(\bfx,\thetavec,t)$ even in the presence of unknown disturbances
$\varepsilon(t)\in L_2^1[t_0,\infty]\cap L_\infty^1[t_0,\infty]$.
These characterizations are the consequence of the fact that
$(f(\bfx(t),\thetavec,t)-f(\bfx(t),\hat{\thetavec}(t),t)))\in
L_{2}^1 [t_0,\infty]$, which in turn is guaranteed by properties
(\ref{eq:assume_alpha}), (\ref{eq:assume_gamma}),
(\ref{eq:assume_alpha_upper}) of the function
$f(\bfx,\thetavec,t)$ in Assumptions \ref{assume:alpha},
\ref{assume:alpha_upper}. Among these properties, estimate
(\ref{eq:assume_alpha_upper}) in Assumption
\ref{assume:alpha_upper} is particulary important for allowing
disturbances (potentially unbounded) from $L_2^1[t_0,\infty]$.
When no disturbances are present
%, or information on the upper bound for
%$\alphavec(\bfx,t)$ is available,
it is possible to show that
P1--P4 hold without involving Assumption \ref{assume:alpha_upper}.

\begin{cor}\label{cor:disturbance} Let system (\ref{system1}),
(\ref{error_model_d}),(\ref{fin_forms_ours_tr1}),
(\ref{fin_forms_ours_tr11}) be given, $\varepsilon(t)= 0$, and
Assumptions \ref{assume:alpha},\ref{assume:explicit_realizability}
hold. Then

P5) norm $\|\thetavec-\hat{\thetavec}(t)\|^{2}_{\Gamma^{-1}}$ is
non-increasing  and properties P1--P4\footnote{In this case,
however, the bound for $\|\psi(\bfx(t),t)\|_{\infty,[t_0,\infty]}$
will be different from the one given by equation
(\ref{eq:psi_gain}) in Theorem \ref{stability_theorem}. Its new
estimate is given by formula (\ref{eq:psi_gain_e_zero}) in
Appendix 1} of Theorem \ref{stability_theorem} hold with
$\varepsilon(t)= 0$ respectively.

%P6) Let $\varepsilon(t)\in L_\infty^1$ and $\alphavec(\bfx,t)$ be
%globally bounded $\|\alphavec(\bfx,t)\|\leq D_\alpha, \
%D_\alpha\in \Real_+$.
%Then
%\begin{equation}\label{eq:disturbance_cor}
%\lim\sup_{T\rightarrow\infty}\frac{1}{T}
%\int_{0}^T(f(\bfx(\tau),\thetavec,\tau)-f(\bfx(\tau),\hat{\thetavec}(\tau),\tau))^2d\tau\leq
%D_\alpha \|\varepsilon(t)\|_{\infty}
%\end{equation}
\end{cor}

In addition to the fact that
$|f(\bfx,\thetavec,t)-f(\bfx,\hat{\thetavec},t)|$ is not required
to be bounded from below as in (\ref{eq:assume_alpha_upper}),
Corollary \ref{cor:disturbance} ensures that
$\|\thetavec-\hat{\thetavec}(t)\|^{2}_{\Gamma^{-1}}$ is {\it not
growing} with time when $\varepsilon(t)=0$.
% It also provides us
%with estimate (\ref{eq:disturbance_cor}) which can further be used
%as a measure of convergence {\it in average}.
The practical relevance of the corollary is that it will allow us
to guarantee desired convergence (\ref{eq:convergence}) with a
much weaker, local version of Assumption \ref{assume:alpha_upper}.
It will also help us to establish conditions for (semi-global)
exponential stability in the unperturbed system, which in turn
will enable (small) disturbances from $L_{\infty}^1[t_0,\infty]$
in the right-hand side of (\ref{error_model_d}).

Another consequence of Theorem \ref{stability_theorem} concerns
the specific case when $\varepsilon(t)\in
L_{\infty}^1[t_0,\infty]\cap L_{2}^1[t_0,\infty]$.

% and function
%$f(\bfx,\thetavec,t)$ is globally bounded in $\bfx$ and $t$.

\begin{cor}\label{cor:gains} Let system (\ref{system1}),
(\ref{error_model_d}),(\ref{fin_forms_ours_tr1}),
(\ref{fin_forms_ours_tr11}) be given, Assumptions
\ref{assume:psi},
\ref{assume:alpha}--\ref{assume:explicit_realizability} hold,
$\varepsilon(t)\in L_{\infty}^1[t_0,\infty]\cap
L_{2}^1[t_0,\infty]$, and property
H\ref{hyp:globally_Lipschitz_uniform} holds. Let in addition,
system (\ref{eq:target_dynamics}) has $L_{p}^1 [t_0,\infty]\mapsto
L_{\infty}^1[t_0,\infty]$, $p\geq 2$ gain. Then

P6) $\psi(\bfx(t),t)\in L_\infty^1[t_0,\infty]$, $\bfx(t)\in
L_\infty^n[t_0,\infty]$;

P7) if properties  H\ref{hyp:locally_bound_uniform_f},
H\ref{hyp:locally_bound_uniform_phi} hold, and system
(\ref{eq:target_dynamics}) has $L_{p}^1 [t_0,\infty]\mapsto
L_{q}^1 [t_0,\infty]$, $q>1$ gain with respect to input $\zeta(t)$
and state $\psi$, then
$\lim_{t\rightarrow\infty}\psi(\bfx(t),t)=0$;

If in addition property H\ref{hyp:locally_bound_uniform_df} is
satisfied, functions $\alphavec(\bfx,t)$, $\pd \psi(\bfx,t)/\pd t$
are locally bounded with respect to $\bfx$ uniformly in $t$, then
limiting relation (\ref{eq:convergence_f_theorem}) holds as well.
\end{cor}
Corollary \ref{cor:gains} extends applicability of algorithms
(\ref{fin_forms_ours_tr1}), (\ref{fin_forms_ours_tr11}) to systems
(\ref{error_model_d}) with defined $L_{p}^{1}[t_0,\infty]\mapsto
L_\infty^1[t_0,\infty]$ gains for arbitrary $p\geq 2$.
%The estimate of the
%upper bound $\|\psi(\bfx(t),t)\|_\infty$ as a function of the
%initial conditions and parameters is provided in Appendix in
%formula (?).

%Theorem \ref{stability_theorem} and Corollaries
%\ref{cor:disturbance}, \ref{cor:gains} provide sufficient
%conditions  ensuring that $\bfx(t)\in L_\infty^n$ in the
%closed-loop system for a variety of the feedbacks (\ref{control}).
%The closed loop system should not necessarily be stable in
%Lyapunov sense. Neither the target dynamics
%(\ref{eq:target_dynamics}) should be stable. It is required
%instead that system (\ref{eq:target_dynamics}) has $L_{p}^1
%[t_0,\infty]\mapsto L_\infty^1$, $p\geq 2$ gain with respect to the
%additive input $\zeta(t)$. In addition, the theorem and
%corollaries provide conditions that guarantee that
%(\ref{eq:convergence_theorem}), (\ref{eq:convergence_psi_theorem})
%hold. These properties will allow us to show convergence of the
%estimates $\hat{\thetavec}(t)$ to $\thetavec$.

Let us formulate conditions ensuring convergence of the estimates
$\hat{\thetavec}(t)$ to $\thetavec$ in the closed loop system
(\ref{system1}), (\ref{error_model_d}),
(\ref{fin_forms_ours_tr1}), (\ref{fin_forms_ours_tr11}). When the
mathematical model of the uncertainties is linear in its
parameters, i.e.
$f(\bfx,\thetavec,t)=\zetavec(\bfx,t)^T\thetavec$, the usual
requirement for convergence is that signal $\zetavec(\bfx(t),t)$
is {\it persistently  exciting} \cite{Sastry89}:
\begin{defn}[Persistent Excitation]\label{defn:LPE} Let function
$\zetavec:\Real_+\rightarrow\Real^k$ be given. Function
$\zetavec(t)$ is said to be persistently exciting iff there exist
constants $\delta>0$ and $L>0$ such that for all $t\in \Real_+$
the following holds
\begin{equation}\label{eq:LPE}
\int_t^{t+L}\zetavec(\tau)\zetavec(\tau)^T d\tau\geq \delta I
\end{equation}
\end{defn}
The conventional notion of persistent excitation requires specific
properties (i.g. the integral inequality (\ref{eq:LPE})) from
signal $\zetavec(t)$ as a function of time. In the closed loop
system, however, relevant signals in the model of uncertainty
$f(\bfx,\thetavec,t)$ can depend on state $\bfx$ and parameters.
In particular, they depend on initial conditions, parameters of
the feedback, and initial time $t_0$. In order to address this
issue it is suggested in \cite{Lorea_2002} to use the notion of
uniform persistent excitation:
\begin{defn}[Uniform Persistent Excitation]\label{defn:LPE_uniform} Let
function  $\zetavec:\Real^n\times\Real_+\rightarrow\Real^k$ be
given, and $\bfx(t,\bfx_0,t_0,\thetavec_0)$ be a solution of
(\ref{system1}), where the vector $\thetavec_0\in\Real^s$ stands
for all possible parameters of  (\ref{system1}) and feedback
(\ref{control}), (\ref{fin_forms_ours_tr1}),
(\ref{fin_forms_ours_tr11}). Function
$\zetavec(\bfx(t,\bfx_0,t_0,\thetavec_0),t)$ is said to be
uniformly persistently exciting iff there exist constants
$\delta>0$ and $L>0$ such that for all $t,t_0\in \Real_+$,
$\bfx_0\in\Real^n$, $\thetavec_0\in\Real^s$ the following holds
\begin{equation}\label{eq:LPE_uniform}
\int_t^{t+L}\zetavec(\bfx(\tau,\bfx_0,t_0,\thetavec_0),\tau)\zetavec(\bfx(\tau,\bfx_0,t_0,\thetavec_0),\tau)^T
d\tau\geq \delta I
\end{equation}
\end{defn}
When dealing with nonlinear parameterization, it is also useful to
have a characterization which takes into account nonlinearity in
the model. In the linear case, persistent excitation of signal
$\zetavec(\bfx(t),t)$ (inequality (\ref{eq:LPE})) implies that the
following property holds
\begin{equation}\label{eq:LPE_1}
\exists \ t'\in[t,t+L]: \
|\zetavec(\bfx(t'),t')^T(\thetavec_1-\thetavec_2)|\geq \delta
\|\thetavec_1-\thetavec_2\|
\end{equation}
In the other words, the difference
$|\zetavec(\bfx(t),t)^T(\thetavec_1-\thetavec_2)|$ is {\it
proportional} to the distance $\|\thetavec_1-\thetavec_2\|$ in
parameter space for some $t'\in[t,t+L]$. In the nonlinear case it
is natural to replace the linear term
$\zetavec(\bfx(t'),t')^T(\thetavec_1-\thetavec_2)$ in
(\ref{eq:LPE_1}) with its nonlinear substitute
$f(\bfx(t'),\thetavec_1,t')-f(\bfx(t'),\thetavec_2,t')$ as has
been done, for example, in  \cite{Cao_2003} for systems with
convex/concave parametrization. It is also natural to replace the
proportion $\delta\|\thetavec_1-\thetavec_2\|$ in the right-hand
side of (\ref{eq:LPE_1}) with a nonlinear function. Therefore, as
a candidate for the {\it nonlinear} persistent excitation
condition we propose the following notion:

\begin{defn}[Nonlinear Persistent Excitation]\label{defn:NlPE} The
function $f(\bfx(t),\thetavec,t):
\Real^n\times\Real^d\times\Real_+\rightarrow\Real$ is said to be
persistently excited with respect to parameters
$\thetavec\in\Omega_\theta\subset\Real^d$ iff there exist constant
$L>0$ and function $\varrho:\Real_+\rightarrow\Real_+, \ \rho\in
\mathcal{K}\cap C^0$ such that for all $t\in \Real_+$,
$\thetavec_1,\thetavec_2\in\Omega_\theta$ the following holds:
\begin{equation}\label{eq:NLPE_2}
\exists \ t'\in[t,t+L]: \
|f(\bfx(t'),\thetavec_1,t')-f(\bfx(t'),\thetavec_2,t')|\geq
\varrho(\|\thetavec_1-\thetavec_2\|)
\end{equation}
\end{defn}

Properties (\ref{eq:LPE}) and (\ref{eq:NLPE_2}) in Definitions
\ref{defn:LPE} and \ref{defn:NlPE} can be considered as
alternative characterizations of excitation in dynamical systems.
While inequality (\ref{eq:LPE}) accounts for specific properties
of the signals  in the uncertainty, inequality (\ref{eq:NLPE_2})
accounts for possibility to detect parametrical difference from
the difference
$f(\bfx(t),\thetavec_1,t)-f(\bfx(t),\thetavec_2,t)$. Taking into
account these two equally possible but still rather distinct
characterizations of excitation in nonlinear systems, in Theorem
\ref{convergence_theorem} below we present a set of alternatives
for parameter convergence in system (\ref{system1}),
(\ref{error_model}), (\ref{fin_forms_ours_tr1}),
(\ref{fin_forms_ours_tr11}).

\begin{thm}[Convergence]\label{convergence_theorem} Let system
(\ref{system1}), (\ref{error_model}), (\ref{fin_forms_ours_tr1}),
(\ref{fin_forms_ours_tr11}) satisfy Assumptions
\ref{assume:psi}--\ref{assume:alpha}. Let, in addition, Assumption
\ref{assume:explicit_realizability} hold with
$\mathcal{B}(\bfx,t)=0$. Then $\bfx(t)\in L^n_\infty[t_0,\infty]$,
$\hat{\thetavec}(t)\in L_\infty^d[t_0,\infty]$. Moreover the
limiting relation:
\[
\lim_{t\rightarrow\infty}\hat{\thetavec}(\bfx(t),t)=\thetavec
\]
is ensured if $\alphavec(\bfx,t)$ is locally bounded in $\bfx$
uniformly in $t$,  and one of the following alternatives hold:

1) function $\alphavec(\bfx(t),t)$ is persistently exciting, and
hypothesis H\ref{hyp:local_linear_bound} holds;

2) function $f(\bfx(t),\thetavec,t)$ is nonlinearly persistently
exciting, i. e. it satisfies condition (\ref{eq:NLPE_2}); it
satisfies hypotheses H\ref{hyp:locally_bound_uniform_f},
H\ref{hyp:locally_bound_uniform_df}; function
$\varphi(\psi,\omegavec,t)$ satisfies
H\ref{hyp:locally_bound_uniform_phi}; function $\pd
\psi(\bfx,t)/\pd t$ be locally bounded in $\bfx$ uniformly in $t$;

In case alternative  1) is satisfied, the estimates
$\hat{\thetavec}(\bfx(t),t)$ converge to $\thetavec$ exponentially
fast. If, in addition, $\alphavec(\bfx(t),t)$ is uniformly
persistently exciting and Assumption \ref{assume:alpha_upper}
holds, then convergence is uniform. The rate of convergence can be
estimated as follows:
\begin{equation}\label{eq:parameters_rate_theorem}
\|\hat{\thetavec}(t)-\thetavec\|\leq  e^{-\rho t}
\|\hat{\thetavec}(t_0)-\thetavec\| D_{\Gamma}
\end{equation}
\[
\rho=\frac{\delta
D_1\lambda_{\min}(\Gamma)}{2L(1+\lambda_{\max}^2(\Gamma)L^2 D^2
\alpha_\infty^4)}, \  \  \ D_{\Gamma}=
\left(\frac{\lambda_{\max}(\Gamma)}{\lambda_{\min}(\Gamma)}\right)^{\frac{1}{2}},
\
\alpha_\infty=\sup_{\|\bfx\|\leq\|\bfx(t)\|_{\infty,[t_0,\infty]},
\ t\geq t_0}\|\alphavec(\bfx,t)\|
\]
\end{thm}

Notice that Theorem \ref{convergence_theorem} considers error
models (\ref{error_model}) where no disturbance term
$\varepsilon(t)$ is present. Despite this Theorem
\ref{convergence_theorem} can be straightforwardly extended to
error models with disturbance (\ref{error_model_d}). Indeed, as
follows from alternative 1), the parameter estimation subsystem
becomes exponentially stable in case function $\alpha(\bfx(t),t)$
is (uniformly) persistently exciting. This in turn allows
(sufficiently small) additive disturbances in the right-hand side
of (\ref{error_model}). In case the excitation is uniform,
convergence of the estimates $\hat{\thetavec}(t)$ to a
neighborhood of $\thetavec$ is guaranteed for every
$\varepsilon(t)\in L_\infty^1[t_0,\infty]$ by inverse Lyapunov
stability theorems \cite{Khalil_2002}.

In case of alternative 2), nonlinear persistent excitation
condition (\ref{eq:NLPE_2}) guarantees convergence
(\ref{eq:convergence}) without invoking Assumption
\ref{assume:alpha_upper} or H\ref{hyp:local_linear_bound}. In this
case, however, the convergence may not be robust, which seems to
be a natural tradeoff between generality of nonlinear
parameterizations $f(\bfx,\thetavec,t)$ and robustness with
respect to unknown disturbances $\varepsilon(t)$.

\section{Discussion}

So far we have shown that, for the class of nonlinearly
parameterized systems, there exist a control function and
parameter adjustment algorithms, such that solutions of the whole
system are bounded and parametric uncertainty is decreasing in
time. We have shown also that in case of persistently excited
functions $\alphavec(\bfx,t)$ the estimates $\hat{\thetavec}(t)$
in (\ref{fin_forms_ours_tr1}) converge exponentially fast to
vector $\thetavec$. These results, however, are not necessarily
limited to  functions satisfying Assumptions \ref{assume:alpha} or
\ref{assume:alpha_upper}.  Due to  space limitations, however, we
provide just the main ideas of possible extensions, leaving out
the technical details. Let us first examine the case where these
assumptions hold only in some domains of the system state space.

{\it Nonlinear functions satisfying Assumptions
\ref{assume:alpha}, \ref{assume:alpha_upper} in a domain of
$\Real^n$}.  Let, in particular, for the given nonlinear function
$f(\bfx,\thetavec,t)$ there exits the following partition of the
state space:
\[
\Omega_{\bfx}=\Omega_M\cup\Omega_{A}, \ \ \Omega_M=\bigcup_{j}
\Omega_{M,j}, \ \Omega_{A}=\Omega_{\bfx}/\Omega_{M}
\]
where $\Omega_{M,j}$ are the  subsets of $\Real^{n}$ where
Assumptions \ref{assume:alpha}, \ref{assume:alpha_upper} are
satisfied for every $\thetavec\in\Omega_\theta$ with the
corresponding functions $\alphavec_j(\bfx,t)$ and constants $D_j$,
$D_{1,j}$. Let us also assume that $\Omega_M$ contains an open
set.

A typical example of a nonlinear function which satisfies this
assumption is $\sin(\theta x)$, where the unknown parameter
$\theta$ belongs to a bounded interval. Let, for instance, the
system dynamics is given by
\begin{equation}\label{ex:sine}
\begin{split}
\dot{x}_1=&x_2\\
\dot{x}_2=&\sin(\theta x_1)+u,
\end{split}
\end{equation}
where parameter $\theta\in \Omega_\theta=[0.6,1.4]$ is unknown
a-priori. For the given bounds of $\Omega_\theta$ the domain
$\Omega_M$ can be derived as follows:
\begin{equation}\label{ex:sine_domain}
\begin{split}
\Omega_M & = \{\bfx \ | \ x_1\in [-3.38,-2.59]\}\cup \{\bfx \ | \
x_1\in [-1.14,1.14]\}\cup \{\bfx \ | \ x_1\in [2.59,3.38]\}\\
& =\Omega_{M,1}\cup \Omega_{M,2}\cup \Omega_{M,3}
\end{split}
\end{equation}
and the function $\alphavec(\bfx,t)$, satisfying Assumptions
\ref{assume:alpha}, \ref{assume:alpha_upper} in $\Omega_M$ is
defined as
\[
\alphavec(\bfx,t)=\left\{\begin{array}{ll}
                          x_1, & \bfx\in \Omega_{M,2}\\
                          -x_1 & \bfx\in \Omega_{M,1}\cup\Omega_{M,3}
                          \end{array}
                        \right.
\]
\begin{figure}\label{fig:domains}
\begin{center}
\includegraphics[width=220pt]{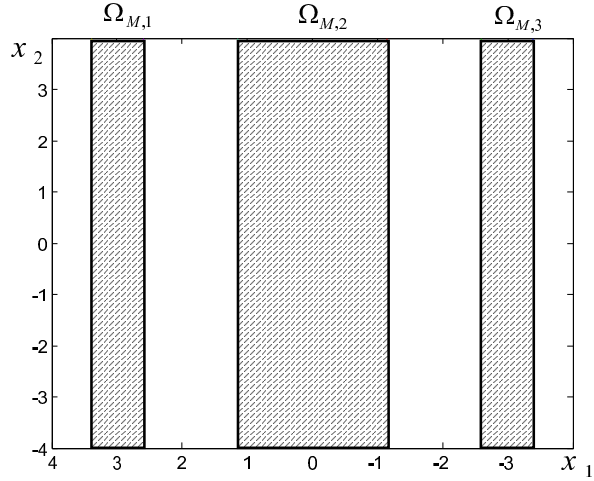}
\end{center}
\caption{Domain $\Omega_M$ for system (\ref{ex:sine}) with
nonlinear parameterized function $\sin(\theta x_1)$,
$\theta\in\Omega_\theta$}
\end{figure}
Another example is $x^{\theta}$, $\theta\in[t_0,\infty)$. The last
parametrization is widely used in modelling physical ``power law"
phenomena in nature (see, for example \cite{Wu_1990}, where this
function models effects of nonlinear damping in muscles).

The fact that Assumptions \ref{assume:alpha},
\ref{assume:alpha_upper} hold in the domain
$\Omega_{M}\subset\Real^n$, allows us to guarantee decrease of the
norm $\|\thetavec-\hat{\thetavec}(t)\|^2_{\Gamma^{-1}}$ only if
the state belongs to $\Omega_M$. Therefore, extra control is
needed in order to steer state $\bfx$ back into the domain
$\Omega_M$. Let us, for example, pick point
$\bfx^\ast\in\Omega_M$, such that $\dist\{\bfx^\ast,\Omega_A\}>r,
\ r\in\Real_+$. Let, in addition, there exists control function
$u_0(\bfx,t)$ such that it steers state $\bfx$ of system
(\ref{system1}) from the initial point $\bfx_0$ into the
$\delta_0$-neighborhood $U(\delta_0,\bfx^\ast)$ of $\bfx^{\ast}$
in finite time $T_0(\bfx_0)$. Suppose also that $\delta_0<r$. As
follows from Theorems \ref{stability_theorem},
\ref{convergence_theorem}, $\hat{\thetavec(t)}$ is bounded for
every segment of the solution which starts from
$U(\delta_0,\bfx^\ast)$ at $t=t_i$ and leaves the domain
$U(r,\bfx^\ast)$ at $t=t_{i+1}$. Furthermore the bound for $
\|\hat{\thetavec}\|_{\infty,[t_i,t_{i+1}]}$ can be estimated
a-priory from the bounds on $\thetavec$ and
$\|\varepsilon(t)\|_{2,[t_0,\infty}$ (see also
(\ref{eq:f_diff_L2})). Given that the right-hand side of system
(\ref{system1}), (\ref{control}), (\ref{fin_forms_ours_tr1}),
(\ref{fin_forms_ours_tr11}) is locally bounded we can conclude
that the time interval $t_{i+1}-t_{i}$ will always be separated
from zero. Taking into account  the results of Theorem
\ref{convergence_theorem}, equation
(\ref{eq:parameters_rate_theorem}), we may conclude that
sufficiently high excitation, defined by the ratio $\delta/L$,
will guarantee that $\|\hat{\thetavec}(t_{i+1})-\thetavec\|<\kappa
\|\hat{\thetavec}(t_i)-\thetavec\|$, $\kappa\in\Real, \
0<\kappa<1$. If the time sequence $\{t_i\}$ is infinite (i.e. the
system always escapes the ball $U(r,\bfx^\ast)$) then convergence
is asymptotic. In case the sequence $\{t_i\}$ is finite (i.e.
$\exists \ t^\ast>0: \ \bfx(t)\in\Omega_M \ \forall t>t^\ast$)
convergence is exponential, this follows from Theorem
\ref{convergence_theorem}.

In principle, the size of $\Omega_M$ and its location in
$\Real^{n}$ depend on the bounds of $\Omega_\theta$. In fact, the
larger the bounds, the smaller the volume of $\Omega_M$. Moreover,
the size of $\Omega_M$ as a function of the bounds of
$\Omega_\theta$ depends on specific properties of nonlinearity
$f(\bfx,\thetavec,t)$. These observations suggest that in order to
handle a broader class of nonlinearities (or functions with higher
degree of uncertainty in $\thetavec$) within the strategy proposed
above, one needs to increase the excitation in functions
$\alphavec_j(\bfx,t)$. This is consistent with previously reported
results \cite{Cao_2003} on parameter convergence in nonlinearly
parameterized systems. Whether the extension of the class of
nonlinearities to more general functions renders it necessary to
increase excitation, however, is still an open issue\footnote{An
example is given in \cite{Cao_2003}, where nonlinear persistent
excitation condition holds for the given parametrization, while
linear persistent excitation condition for linear parametrization
with respect to the same parameter-independent function is not
satisfied.}.

{\it Functions $f(\bfx,\thetavec,t)$ with nonlinear incremental
growth rates in $\thetavec$}. Another direction to extend the
class nonlinear functions suitable for our method is to allow
nonlinear bounds for the growth rates in (\ref{eq:assume_alpha}),
(\ref{eq:assume_alpha_upper}) in Assumptions \ref{assume:alpha},
\ref{assume:alpha_upper}. The most straightforward
generalizations, which do not change dramatically the machinery of
technical proofs of Theorems \ref{stability_theorem},
\ref{convergence_theorem}, are provided in Assumptions
{\ref{assume:alpha1}}, \ref{assume:alpha_upper1} below:
\begin{assume}\label{assume:alpha1}For the given function $f(\bfx,\thetavec,t)$ in (\ref{error_model_d}) there exist function
$\alphavec(\bfx,t): \Real^{n}\times \Real_+\rightarrow \Real^d, \
\alphavec(\bfx,t)\in \mathcal{C}^1$, function
$\sigmavec:\Real^d\rightarrow\Real^d$,
$\sigmavec(\bfz)=(\sigma_1(z_1),\sigma_2(z_2),\dots,\sigma_d(z_d))^T$
\[
\sigma_i(\xi)\xi \geq 0, \
\lim_{S\rightarrow\infty}\int_0^{S}\sigma_i(\xi)d\xi=\infty,
\]
and function $\overline{\gamma}\in\mathcal{K}$ such that
\begin{equation}\label{eq:assume_alpha1}
(f(\bfx,\hat{\thetavec},t)-f(\bfx,\thetavec,t))\alphavec(\bfx,t)^{T}\sigmavec(\hat{\thetavec}-\thetavec)\geq0
\end{equation}
\begin{equation}\label{eq:assume_gamma1}
|f(\bfx,\hat{\thetavec},t)-f(\bfx,\thetavec,t)|\leq
\overline\gamma(|\alphavec(\bfx,t)^T\sigmavec(\thetavec-\hat{\thetavec})|)
\end{equation}
\end{assume}

\begin{assume}\label{assume:alpha_upper1} For the given function $f(\bfx,\thetavec,t)$ in
(\ref{error_model_d}) and function $\alphavec(\bfx,t)$ satisfying
Assumption \ref{assume:alpha1} there exists function
$\underline{\gamma}\in\mathcal{K}$ such that
\[
\underline\gamma(|\alphavec(\bfx,t)^T\sigmavec(\thetavec-\hat{\thetavec})|)\leq|f(\bfx,\thetavec,t)-f(\bfx,\hat{\thetavec},t)|
\]
\end{assume}

Choosing for simplicity $\Gamma=I$,  denoting $\Delta
f(\hat{\thetavec},\thetavec,\bfx,t)=
|f(\bfx,\hat{\thetavec},t)-f(\bfx,\thetavec,t)|$, letting
$|\varepsilon(t)|$ in (\ref{error_model_d}) be such that
$\int_{0}^{\infty}\gamma_{\epsilon}(|\varepsilon(\tau)|)d\tau<\infty$,
$\gamma_\epsilon\in\mathcal{K}$ and replacing
$V_{\hat{\thetavec}}(\hat{\thetavec},\thetavec,t)$ in
(\ref{V_theta}) (proof, of Theorem \ref{stability_theorem},
Appendix 1) with
\[
V_{\hat{\thetavec}}(\hat{\thetavec},\thetavec,t)=\sum_{i=1}^d\int_0^{\sigma(\hat{\theta}_i-\theta_i)}\sigma(\xi)d\xi+\int_t^\infty\gamma_{\epsilon}(|\varepsilon(\tau)|)d\tau
\]
we can obtain the following estimate:
\begin{eqnarray}\label{eq:nonlinear_gains_1}
\dot{V}&=&-\sigmavec(\hat{\thetavec}-\thetavec)^T\alphavec(\bfx,t)(f(\bfx,\hat{\thetavec},t)-f(\bfx,\thetavec,t))+\varepsilon(t)\sigmavec(\hat{\thetavec}-\thetavec)^T\alphavec(\bfx,t)-\gamma_\epsilon(|\varepsilon(t)|)\nonumber \\
&\leq& -\overline{\gamma}^{-1}(\Delta
f(\hat{\thetavec},\thetavec,\bfx,t))\Delta
f(\hat{\thetavec},\thetavec,\bfx,t)+|\varepsilon(t)|\underline{\gamma}^{-1}(\Delta
f(\hat{\thetavec},\thetavec,\bfx,t))-\gamma_\epsilon(|\varepsilon(t)|)
\end{eqnarray}
Boundedness of $\hat{\thetavec}$ follows from
(\ref{eq:nonlinear_gains_1}) if we can resolve the following
inequality for unknown $\gamma_\epsilon(|\varepsilon(t)|)$:
\[
-\overline{\gamma}^{-1}(\Delta
f(\hat{\thetavec},\thetavec,\bfx,t))\Delta
f(\hat{\thetavec},\thetavec,\bfx,t)+|\varepsilon(t)|\underline{\gamma}^{-1}(\Delta
f(\hat{\thetavec},\thetavec,\bfx,t))-\gamma_\epsilon(|\varepsilon(t)|)\leq
0,
\]
If in addition there exists $\gamma_f\in\mathcal{K}$ such that
\begin{equation}\label{eq:nonlinear_gains_2}
-\overline{\gamma}^{-1}(\Delta
f(\hat{\thetavec},\thetavec,\bfx,t))\Delta
f(\hat{\thetavec},\thetavec,\bfx,t)+|\varepsilon(t)|\underline{\gamma}^{-1}(\Delta
f(\hat{\thetavec},\thetavec,\bfx,t))-\gamma_\epsilon(|\varepsilon(t)|)\leq
-\gamma_f(|\Delta f(\hat{\thetavec},\thetavec,\bfx,t)|)
\end{equation}
then we can guarantee that $\Delta
f(\hat{\thetavec}(t),\thetavec,\bfx(t),t)\in
L_{\gamma_f}^1[t_0,\infty]$, where $L_{\gamma_f}^n[t_0,\infty]$ is
the space of all functions $\bff_0(t):\Real_+\rightarrow\Real^n$
with finite integral $\int_{0}^{\infty}\gamma_f(\|\bff_0(t)\|)d t
\leq \infty$. Therefore, if system (\ref{eq:target_dynamics}) has
$L_{\gamma_f}^1[t_0,\infty]\cup L_{\gamma_\epsilon}^1\mapsto
L_\infty^1$ gain, we can conclude, invoking this new modified
Assumption \ref{assume:gain}, that $\psi(\bfx(t),t)\in
L_\infty^1$, $\bfx(t)\in L_\infty^n$.

Notice that letting functions $\underline{\gamma}$,
$\overline{\gamma}$ linear ($\underline{\gamma}(|s|)=D |s|$,
$\overline{\gamma}(|s|)=D_1|s|$), allows us straightforwardly
obtain, as in (\ref{parameric_deviation_derivative}), that choice
$\gamma_\epsilon=\frac{D}{4D_1^2}\varepsilon^2$ ensures the
following inequality:
\[
\dot{V}\leq
-\frac{1}{D}\left(|f(\bfx,\hat{\thetavec},t)-f(\bfx,\thetavec,t)|-\frac{D}{2
D_1} \varepsilon(t)\right)^2
\]
This implies that properties similar to P1)--P7) can be derived
for the case where Assumptions \ref{assume:alpha},
\ref{assume:alpha_upper} are replaced with Assumptions
\ref{assume:alpha1}, \ref{assume:alpha_upper1} and functions
$\overline{\gamma}(\cdot)$, $\underline{\gamma}(\cdot)$  are
linear (for the proofs of Theorem \ref{stability_theorem}, and
Corollaries \ref{cor:disturbance}, \ref{cor:gains} see the
Appendix). Parameter convergence in this case can also be deduced
from Theorem \ref{convergence_theorem}, alternative 2). Notice
that, due to the nonlinearities in
$\sigmavec(\thetavec-\hat{\thetavec})$, convergence in general may
not be exponentially fast.

For the nonlinear functions $\underline{\gamma}$,
$\overline{\gamma}$ in Assumptions \ref{assume:alpha1},
\ref{assume:alpha_upper1} the formulation of the results will
require the notions of $L_{\gamma_f}$ spaces introduced above. The
machinery behind the statements, however, will remain the same.
The practical relevance of these results with nonlinear functions
$\underline{\gamma}$, $\overline{\gamma}$ is that they enable us
to take into account the specific properties of the signal
$\varepsilon(t)$ when designing control, adaptation and parameter
estimation procedures. If, for example, disturbance
$\varepsilon(t)$ is due to observers, we might derive requirements
on convergence rates for the observer-induced errors
$\varepsilon(t)$ (i. e. $\varepsilon(t)\in
L_{\gamma_\epsilon}^1[t_0,\infty]$, and $\gamma_\epsilon(\cdot)$
satisfies inequality (\ref{eq:nonlinear_gains_2})). Given these
rates and the fact that $\Delta f(\bfx(t),\hat{\thetavec}(t),t)\in
L_{\gamma_f}^1 [t_0,\infty]$, the target dynamics
\[
\dpsi=-\varphi(\psi,\omegavec,t)+\zeta(t), \ \zeta(t)\in
L_{\gamma_f}^1 [t_0,\infty]\cup L_{\gamma_\epsilon}^1[t_0,\infty]
\]
should be chosen in order to guarantee boundedness of  $\psi(t)$
for all $\zeta(t)\in L_{\gamma_f}^1 [t_0,\infty]\cup
L_{\gamma_\epsilon}^1[t_0,\infty]$. This will allow synergy at all
stages of the design and analysis of adapting systems.

{\it Singularities in control (\ref{control}), and non-affine
models}. In the problem statement we restricted the class of
nonlinear systems of interest models (\ref{system1}) that are
affine in control and furthermore, we assumed that inverse
$\left(L_{\bfg(\bfx)}\psi(\bfx,t)\right)^{-1}$ exists everywhere.
Even though this restriction holds in wide variety of practically
relevant situations, the question is whether the proposed approach
could be extended to more general classes of systems. Let us, for
instance, assume that either $L_{\bfg(\bfx)}\psi(\bfx,t)=0$ for
some $\bfx\in\Real^n$, or the right-hand side of (\ref{system1})
is not affine in control, e.g.
\begin{equation}\label{eq:non_affine}
\dot{\bfx}=\bff(\bfx,\thetavec,u)
\end{equation}
Obviously, control function (\ref{control}), which transforms
(\ref{system1}) into (\ref{error_model}), is not relevant any
more.  Despite that, it is still possible to transform system
(\ref{system1}) into an error model, similar to
(\ref{error_model}). In order realize this transformation without
invoking the use of linearity in the control or taking inverse
$\left(L_{\bfg(\bfx)}\psi(\bfx,t)\right)^{-1}$, it should be
possible to find a function $u(\bfx,\thetavec,\omegavec,t)$ such
that the following {\it invariance} condition is satisfied:
\begin{equation}\label{eq:invariance_control}
\frac{\pd \psi(\bfx,t)}{\pd \bfx
}\bff(\bfx,\thetavec,u(\bfx,\thetavec,\omegavec,t))=-\varphi(\psi(\bfx,t),\omegavec,t)-\frac{\pd
\psi(\bfx,t)}{\pd t}
\end{equation}
Denoting
\begin{equation}\label{eq:new_f}
\frac{\pd \psi(\bfx,t)}{\pd \bfx
}\bff(\bfx,\thetavec,u(\bfx,\hat{\thetavec},\omegavec,t))=-f^\ast(\bfx,\thetavec,\hat{\thetavec},\omegavec,t)
\end{equation}
and taking into account (\ref{eq:non_affine}),
(\ref{eq:invariance_control}), and (\ref{eq:new_f}) we can
calculate derivative $\dot{\psi}$ in the following form
\begin{eqnarray}\label{eq:general_error_model}
\dot{\psi}&=&-f^\ast(\bfx,\thetavec,\hat{\thetavec},\omegavec,t)+f^\ast(\bfx,\thetavec,{\thetavec},\omegavec,t)-f^\ast(\bfx,\thetavec,{\thetavec},\omegavec,t)+\frac{\pd
\psi(\bfx,t) }{\pd t}\nonumber \\
&=&
-\varphi(\psi(\bfx,t),\omegavec,t)+f^\ast(\bfx,\thetavec,{\thetavec},\omegavec,t)-f^\ast(\bfx,\thetavec,\hat{\thetavec},\omegavec,t)
\end{eqnarray}
The main difference between error models
(\ref{eq:general_error_model}) and (\ref{error_model}) is that
function $f^\ast(\bfx,\thetavec,\hat{\thetavec},\omegavec,t)$ in
(\ref{eq:general_error_model}) depends on additional parameters
$\thetavec$, $\omegavec$. Despite this difference our approach can
still be applied to models (\ref{eq:general_error_model}) if
inequalities (\ref{eq:assume_alpha}),
(\ref{eq:assume_alpha_upper}) in Assumption \ref{assume:alpha}
(or/and \ref{assume:alpha_upper}) hold for function
$f^\ast(\bfx,\thetavec,\hat{\thetavec},\omegavec,t)$ for any
$\omegavec\in\Omega_\omega$. Adaptation algorithms in this case
can straightforwardly be derived from (\ref{t2_1}), (\ref{t2_2})
(in Appendix 1) and will have the form similar to
(\ref{fin_forms_ours_tr1}), (\ref{fin_forms_ours_tr11}).

In the next section we illustrate the application to and main
steps in the design of our algorithms for the optimal slip
identification problem in brake control systems.

\section{Example}

%\subsection{Monotonic parameterized models}

Consider the problem of minimizing  the braking distance for a
single wheel rolling along a surface. The surface properties can
vary depending on the current position of the wheel. The wheel
dynamics can be given by the following system of differential
equations \cite{I_Petersen_2003}:
\begin{eqnarray}\label{ex:slip_model}
\dot{x}_1&=&-\frac{1}{m}F_s(F_n,\bfx,\theta), \nonumber\\
\dot{x}_2&=&\frac{1}{J}(F_s(F_n,\bfx,\theta)r-u)\nonumber \\
\dot{x}_3&=&-\frac{1}{x_1}((\frac{1}{m}(1-x_3)+\frac{r^2}{J})F_s(F_n,\bfx,\theta)-\frac{r}{J}u),
\end{eqnarray}
$x_1$ is longitudinal velocity, $x_2$ is angular velocity,
$x_3=(x_1-r x_2)/x_1$\footnote{Given this functional relation,
variable $x_3$ in (\ref{ex:slip_model}) can be viewed as an output
of the reduced system with state $(x_1,x_2)$.} is wheel slip,
$\bfx=(x_1,x_2,x_3)^T$, $m$ is mass of the wheel, $J$ is moment of
inertia, $r$ is radius of the wheel, $u$ is control input (brake
torque), $F_s(F_n,\bfx,\theta)$ is a function specifying the
tyre-road friction force depending on the surface-dependent
parameter $\theta$ and the load force $F_n$. This function, for
example, can be derived from steady-state behavior of the LuGre
tyre-road friction model \cite{Canudas_1999}:
\begin{equation}\label{eq:ex_friction}
F_s(F_n,\bfx,\theta)=F_n \sign(x_2)\frac{\frac{\sigma_0}{L}
g(x_2,x_3,\theta)\frac{x_3}{1-x_3}}{\frac{\sigma_0}{L}
\frac{x_3}{1-x_3}+g(x_2,x_3,\theta)},
\end{equation}
\begin{equation}\label{eq:ex_friction_component}
g(x_2,x_3,\theta)=\theta(\mu_C + (\mu_{S}-\mu_{C})e^{-\frac{|r x_2
x_3|}{|1-x_3|v_s}}), \ \
\end{equation}
where $\mu_C$, $\mu_S$ are Coulomb and static  friction
coefficients, $v_s$ is the Stribeck velocity, $\sigma_0$ is the
normalized rubber longitudinal stiffness, $L$ is the length of the
road contact patch. In order to avoid singularities we assume, as
suggested in \cite{I_Petersen_2003}, that the system is turned off
when velocity $x_1$ reaches a small neighborhood of zero (in our
example we stopped simulations as soon as $x_1$ becomes less than
$5$ m/sec). Moreover, given that functions (\ref{eq:ex_friction}),
(\ref{eq:ex_friction_component}) are bounded for the relevant set
of the system parameters, it is always possible to design control
function $u(\bfx,t)$ in (\ref{ex:slip_model}) such that
\begin{equation}\label{eq:ex_x3_bounds}
x_3(t)\in [\delta,1-\delta], \ \delta\in\Real, \ \delta>0
\end{equation}
for all $t: \ x_1(t)\geq \delta_1$, $\delta_1=5$ m/sec.

While the majority of the model parameters can be estimated
a-priori, the tyre-road parameter $\theta$ is dependent on the
properties of the road surface. Therefore, on-line identification
of the parameter $\theta$ is desirable in order to compute the
optimal slip value
\begin{eqnarray}\label{ex:slip_value}
x_3^{\ast}=\arg\max_{x_3}{F_s(F_n,\bfx,\theta)}
\end{eqnarray}
which ensures the maximum deceleration force and therefore results
in the shortest braking distance.

The main loop controller is derived in accordance with the
standard certainty-equivalence principle and can be written as
follows:
\[
u(\bfx,\hat{\theta},x_3^\ast)=\frac{J}{r}((\frac{1}{m}(1-x_3)+\frac{r^2}{J})F_s(F_n,\bfx,\hat{\theta})-K_s
x_1(x_3-x_3^\ast)), \ K_s>0
\]
In order to estimate parameter $\theta$ by measuring the values of
variables $x_1,x_2$ and $x_3$, we construct the following
subsystem:
\[
\dot{\hat{x}}_3=-\frac{1}{x_1}((\frac{1}{m}(1-x_3)+\frac{r^2}{J})F_s(F_n,\bfx,\hat{\theta})-\frac{r}{J}u)+(x_3-\hat{x}_3)
\]
and consider the dynamics of error function
$\psi(\bfx,t)=\psi(x_3,\hat{x}_3)=x_3-\hat{x}_3$:
\begin{eqnarray}\label{ex:slip_error_model}
\dot{\psi}=-\psi+\frac{1}{x_1}(\frac{1}{m}(1-x_3)+\frac{r^2}{J})(F_s(F_n,\bfx,\theta)-F_s(F_n,\bfx,\hat{\theta}))
\end{eqnarray}
The desired dynamics of system (\ref{ex:slip_error_model}) is
\begin{equation}\label{eq:ex_desired_dynamics}
\dot{\psi}=-\psi + \xi(t)
\end{equation}
where $\xi(t)$ is to be from $L_{2}^1[t_0,\infty]$.  Let us check
Assumptions \ref{assume:psi}, \ref{assume:gain} for the function
$\psi(\bfx,t)=x_3-\hat{x}_3(t)$ and system
(\ref{eq:ex_desired_dynamics}). Notice first that  state $\bfx$ of
system (\ref{ex:slip_model}) is bounded according to the physical
laws governing the dynamics of (\ref{ex:slip_model}). In addition,
boundedness of $\psi(\bfx,t)$ implies that $\hat{x}_3(t)$ is
bounded. Hence Assumption \ref{assume:psi} holds. System
(\ref{eq:ex_desired_dynamics}), obviously, has
$L_{2}^1[t_0,\infty]\mapsto L_\infty[t_0,\infty]$ gain. We can
conclude that Assumption \ref{assume:gain} also holds. Let us
check Assumptions \ref{assume:alpha}, \ref{assume:alpha_upper}.

Taking into account (\ref{eq:ex_x3_bounds}) we can conclude that
function $\frac{1}{x_1}(\frac{1}{m}(1-x_3)+\frac{r^2}{J})$ in
(\ref{ex:slip_error_model}) is positive and, furthermore, is
separated from zero for all $x_1>\delta_1$. Therefore,  taking
this into account equations (\ref{eq:ex_friction}),
(\ref{eq:ex_friction_component}) we can conclude that function
$\frac{1}{x_1}(\frac{1}{m}(1-x_3)+\frac{r^2}{J})F_s(F_n,\bfx,\theta)$
in (\ref{ex:slip_error_model}) satisfies Assumptions
\ref{assume:alpha}, \ref{assume:alpha_upper} with
\[
\alpha(\bfx,t)=\const=1, \ \forall \ x_1: \  \delta_1<x_1<x_1(t_0)
\]
Therefore, in order to design an estimation scheme satisfying
assumptions of Theorem \ref{convergence_theorem} we shall find
functions $\Psi(\bfx,t)$, $\mathcal{B}(\bfx,t)$ such that
Assumption \ref{assume:explicit_realizability} holds. It is easy
to see that this assumption is satisfied with
$\Psi(\bfx,t)=\const$, and $\mathcal{B}(\bfx,t)$=0. Let us choose,
therefore, $\Psi(\bfx,t)=0$.  Then according to
(\ref{fin_forms_ours_tr1}) and (\ref{ex:slip_error_model}), a
parameter adjustment algorithm will be given by the following
system:
\begin{eqnarray}\label{ex:slip_alg}
\hat{\theta}=\gamma((x_3-\hat{x}_3)+\hat{\theta}_I), \
\dot{\hat{\theta}}_I=(x_3-\hat{x}_3), \ \gamma=100
\end{eqnarray}
An important fact about algorithm (\ref{ex:slip_alg}) is that it
is a {parametric} {\it linear proportional-integral} scheme.
According to Theorem \ref{convergence_theorem} the estimates
(\ref{ex:slip_alg}) converge to $\theta$ exponentially fast in the
domain specified by equation (\ref{eq:ex_x3_bounds}), and
inequality $x_1(t_0)\geq x_1(t)>\delta_1$. The last inequality is
satisfied as, according to (\ref{ex:slip_model}), time-derivative
of the variable $x_1(t)$ is non-positive and the system is turned
``off`` when $x_1(t)\leq \delta_1$.

We simulated system (\ref{ex:slip_model}) -- (\ref{ex:slip_alg})
with the following setup of parameters and initial conditions:
$\sigma_0=200$, $L=0.25$, $\mu_C=0.5$, $\mu_S=0.9$, $v_s=12.5$,
$r=0.3$, $m=200$, $J=0.23$, $F_n=3000$, $K_s=30$. The
effectiveness of estimation algorithm (\ref{ex:slip_alg}) could be
illustrated with Figure 4. Estimates $\hat{\theta}$ approach the
actual values of parameter $\theta$ sufficiently fast for the
controller to calculate the optimal slip value $x_3^\ast$ and
steer the system toward this point in real braking time.
Effectiveness of the proposed identification-based control can be
confirmed by comparing the braking distance in the system with
on-line estimation of $x_3^\ast$ according to formula
(\ref{ex:slip_value})  with the one, in which the values of
$x_3^\ast$ were kept constant (in the interval $[0.1,0.2]$). For
model parameters as presently given and road condition given by
the piece-wise constant function
\[
\theta(s)=\left\{
                    \begin{array}{ll}
                     0.3, & s\in[0,8]\\
                     1.3, & s\in(8,16]\\
                     0.7, & s\in(16,24]\\
                     0.4, & s\in(24,32]\\
                     1.5, & s\in(32,40]\\
                     0.6, & s\in(40,\infty]
                    \end{array}
                  \right., s=\int_0^{t}x_1(\tau)d\tau
\]
the simulated braking distance obtained with our on-line
estimation procedure of $x_3^\ast$ is $54.95$ meters. This result
compares favorably with the values obtained for preset values of
$x_3^\ast$, which range between $57.52$ and $55.32$ (for
$x_3^\ast=0.1$ and $x_3^\ast=0.2$ respectively).

\begin{figure}\label{fig:slip}
\begin{center}
\includegraphics[width=200pt]{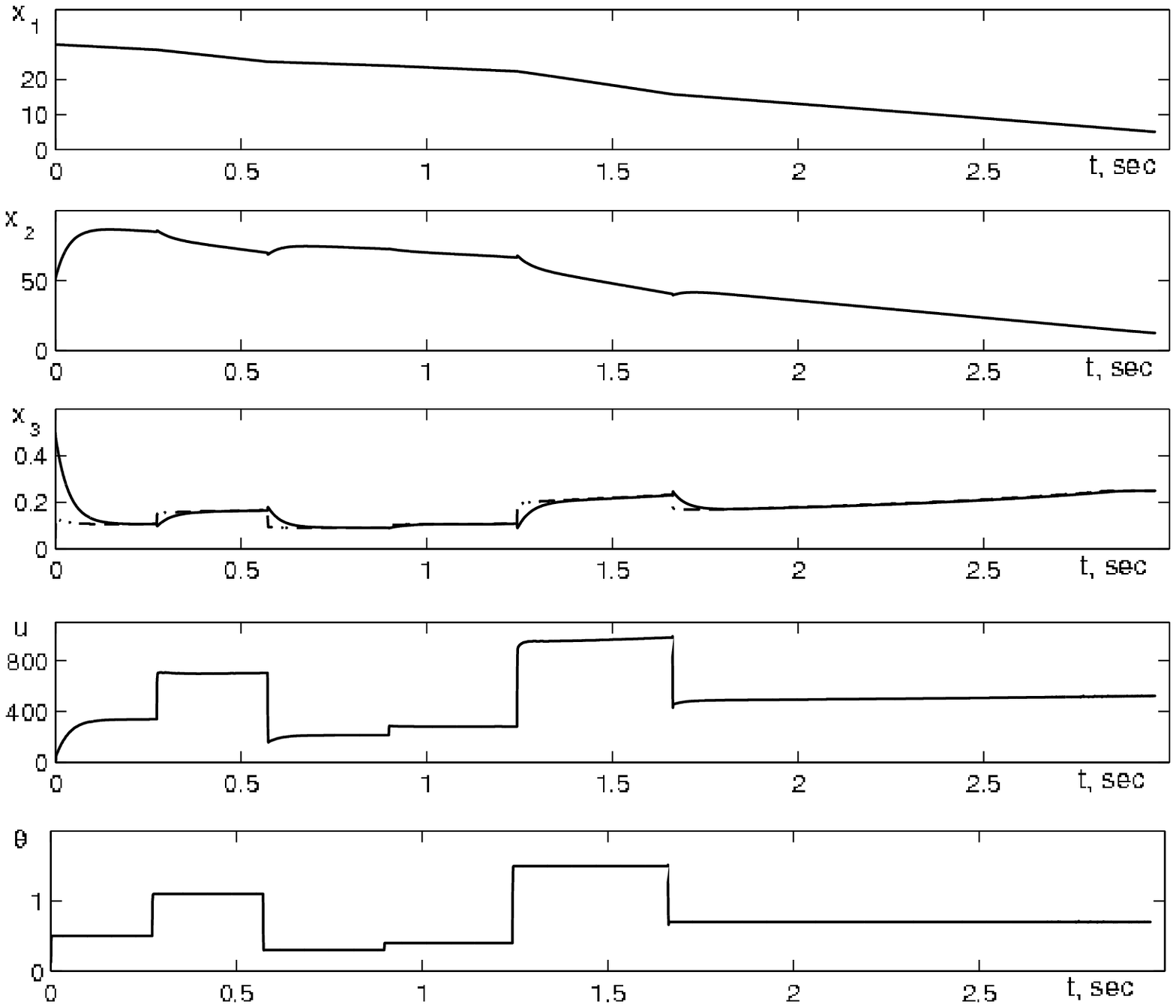}
\end{center}
\caption{Plots of the trajectories of system
(\ref{ex:slip_model})}
\end{figure}

\section{Conclusion}

In the present article we provided new tool for the design and
analysis of adaptive/parameter estimation schemes for dynamic
systems with possible Lyapunov-unstable desired dynamics and
nonlinear parameterization. In our method we consider adaptation
as a {\it process of asymptotic compensation} of the uncertainty,
or as {\it control in functional spaces}, rather than as simply
reaching of a control goal. In particular, we wished to achieve
that mismatches between the modeled uncertainty and compensator
vanish asymptotically with time or { belong} to  specific {
functional spaces}. This understanding of adaptation naturally
leads to the possibility to describe the { desired dynamics} of
adapting systems in terms of an {\it operator}, which maps these
mismatches into  error functions, as functions of time from a
functional space. Continuity of this {target operator} is not
required. Hence stability of the desired dynamics, as a substitute
of continuity, is not necessary for our approach. The adaptation
mechanism itself could be viewed as control in  functional spaces.
In the other words, the aim of adaptation consists in ensuring
that the uncertainty-induced errors of the compensator belong to a
specific functional space. This  idea leads to classes of adaptive
systems, where  applications require gentle, non-dominating
control and where the desired dynamical state can be unstable.

When the desired motions in the system are known to be  Lyapunov
stable, our approach allows to design adaptation procedures
without knowledge of the particular Lyapunov function. As
mentioned in \cite{Fradkov99}, this was one of the open
theoretical challenges in the theory of adaptive cotnrol.

Another contribution of our present study is that we proposed a
new class of parameterizations for nonlinearly parameterized
models. Instead of aiming at a general solution for the problem of
nonlinearity in the parameters, parametrization was restricted to
a set of smooth functions, which are monotonic with respect to a
linear functional in the parameters. For this new class,
adaptation/estimation algorithms were introduced and analyzed. It
was shown  that standard linear persistent excitation conditions
suffice to ensure exponentially fast convergence of the estimates
to the actual values of unknown parameters. If, however, the
monotonicity assumption holds only locally in the system state
space, excitation with sufficiently high-frequency of oscillations
still is able to ensure cpnvergence. In addition to the analysis
of the effects of conventional persistent excitation on
convergence, we also formulated a much weaker property - {\it
nonlinear persistent excitation} condition. With this new property
we established conditions for {\it asymptotic} convergence of the
estimates. It is also desirable to notice that in case of linear
parametrization the proposed parameter estimation schemes  allow
to estimate the unknowns in a dynamical system without asking for
the usual filtered transformations, thus reducing the number of
integrators in the estimator.

An application of our results, which is relevant to the parameter
estimation problems for systems with nonlinear parameterization,
was provided as an example. In this example we did not cover all
solutions to every theoretical problem we were targeting in this
article. In particular, it covers only the problem of nonlinear
parameterization. The main rationale, however, was to illustrate
all steps of our  method. Last but not the least, the application
presents a practically relevant solution to an important
engineering problem. The effectiveness of the solution to this
problem leads us to expect that our newly proposed method can
successfully be implemented in a variety of other applications.

\section{Appendix 1. Proofs of the theorems and auxiliary results}

{\it Proof of Theorem \ref{stability_theorem}.} Let us first show
that property P1) holds. Consider solutions of system
(\ref{system1}), (\ref{error_model_d}),
(\ref{fin_forms_ours_tr1}), (\ref{fin_forms_ours_tr11}) passing
through the point $\bfx(t_0)$, $\hat{\thetavec}_I(t_0)$ for
$t\in[t_0,T^\ast]$\footnote{According to the theorem formulation,
interval $[t_0,T^\ast]$ is the interval of existence of the
solutions}. Let us  calculate formally the time-derivative  of
function $\hat{\thetavec}(\bfx,t)$:
$\dot{\hat{\thetavec}}(\bfx,t)=\Gamma({\dot{\hat{\thetavec}}_{P}}+\dot{\hat\thetavec}_I)=\Gamma(\dpsi\alphavec(\bfx,t)+\psi\dot{\alphavec}(\bfx,t)-\dot{\Psi}(\bfx,t)+\dot{\hat\thetavec}_I)$.
Notice that
\begin{eqnarray}\label{t2_1}
&&\psi\dot{\alphavec}(\bfx,t)-\dot{\Psi}(\bfx,t)+\dot{\hat{\thetavec}}_I=\psi(\bfx,t)\frac{\pd
\alphavec(\bfx,t)}{\pd \bfx_1}\dot{\bfx}_1+\psi(\bfx,t)\frac{\pd
\alphavec(\bfx)}{\pd \bfx_2}\dot{\bfx}_2 + \psi(\bfx,t)\frac{\pd
\alphavec(\bfx,t)}{\pd t}-\nonumber\\
& & \frac{\pd \Psi(\bfx,t)}{\pd \bfx_1}\dot{\bfx}_1-\frac{\pd
\Psi(\bfx,t)}{\pd \bfx_2}\dot{\bfx}_2-\frac{\pd \Psi(\bfx,t)}{\pd
t}+\dot{\hat\thetavec}_I
\end{eqnarray}
According to Assumption \ref{assume:explicit_realizability},
$\frac{\pd \Psi(\bfx,t)}{\pd \bfx_2}=\psi(\bfx,t)\frac{\pd
\alphavec(\bfx,t)}{\pd \bfx_2}+\mathcal{B}(\bfx,t)$. Then taking
into account (\ref{t2_1}), we can obtain
\begin{eqnarray}\label{t2_2}
& &
\psi\dot{\alphavec}(\bfx,t)-\dot{\Psi}(\bfx,t)+\dot{\hat{\thetavec}}_I=\left(\psi(\bfx,t)\frac{\pd
\alphavec(\bfx,t)}{\pd \bfx_1}-\frac{\pd \Psi}{\pd \bfx_1
}\right)\dot{\bfx}_1+\psi(\bfx,t)\frac{\pd \alphavec(\bfx,t)}{\pd
t}-\frac{\Psi(\bfx,t)}{\pd t}-\nonumber
\\
& &
\mathcal{B}(\bfx,t)(\bff_2(\bfx,\thetavec)+\bfg_2(\bfx)u)+\dot{\hat{\thetavec}}_I
\end{eqnarray}
Notice that according to the proposed notation we can rewrite the
term $\left(\psi(\bfx,t)\frac{\pd \alphavec(\bfx,t)}{\pd
\bfx_1}-\frac{\pd \Psi}{\pd \bfx_1 }\right)\dot{\bfx}_1$ in the
following form: $\left(\psi(\bfx,t)L_{\bff_1}
\alphavec(\bfx,t)-L_{\bff_1} \Psi(\bfx,t)\right)+
\left(\psi(\bfx,t)L_{\bfg_1} \alphavec(\bfx,t)-L_{\bfg_1}
\Psi(\bfx,t)\right)u(\bfx,\hat{\thetavec},t)$. Hence it follows
from (\ref{fin_forms_ours_tr1}) and (\ref{t2_2}) that
$\psi\dot{\alphavec}(\bfx,t)-\dot{\Psi}(\bfx,t)+\dot{\hat{\thetavec}}_I=\varphi(\psi)\alphavec(\bfx,t)-\mathcal{B}(\bfx,t)(\bff_2(\bfx,\thetavec)-\bff_2(\bfx,\hat{\thetavec}))$.
Therefore derivative $\dot{\hat\thetavec}(\bfx,t)$ can be written
in the following way:
\begin{eqnarray}\label{algorithm_dpsi}
\dot{\hat{\thetavec}}=\Gamma((\dpsi+\varphi(\psi))\alphavec(\bfx,t)-\mathcal{B}(\bfx,t)(\bff_2(\bfx,\thetavec)-\bff_2(\bfx,\hat{\thetavec})))
\end{eqnarray}
Consider the following positive-definite function:
\begin{equation}\label{V_theta}
V_{\hat{\thetavec}}(\hat{\thetavec},\thetavec,t)=
\frac{1}{2}\|\hat{\thetavec}-\thetavec\|^2_{\Gamma^{-1}} +
\frac{D}{4 D_1^2} \int_{t}^\infty\varepsilon^2(\tau)d\tau
\end{equation}
Its time-derivative according to equations (\ref{algorithm_dpsi})
can be obtained as follows:
\begin{eqnarray}\label{eq:dV_full_alg}
\dot{V}_{\hat{\thetavec}}(\hat{\thetavec},\thetavec,t)=(\varphi(\psi)+\dpsi)(\hat{\thetavec}-\thetavec)^{T}\alphavec(\bfx,t)-(\hat{\thetavec}-{\thetavec})^T\mathcal{B}(\bfx,t)(\bff_2(\bfx,\thetavec)-\bff_2(\bfx,\hat{\thetavec}))-
\frac{D}{4 D_1^2}\varepsilon^2(t)
\end{eqnarray}
Let $\mathcal{B}(\bfx,t)\neq 0$, then  consider the following
difference $\bff_2(\bfx,\thetavec)-\bff_2(\bfx,\hat{\thetavec})$.
Applying Hadamard's lemma we represent this difference in the
following way:
\[
\bff_2(\bfx,\thetavec)-\bff_2(\bfx,\hat{\thetavec})= \int_0^1
\frac{\pd \bff_2(\bfx, \bfs(\lambda))}{\pd \bfs} d\lambda
(\thetavec-\hat{\thetavec}), \ \
\bfs(\lambda)=\thetavec\lambda+\hat{\thetavec}(1-\lambda)
\]
Therefore, according to Assumption
\ref{assume:explicit_realizability} function
$(\hat{\thetavec}-{\thetavec})^T\mathcal{B}(\bfx,t)(\bff_2(\bfx,\thetavec)-\bff_2(\bfx,\hat{\thetavec}))$
is positive semi-definite, hence using Assumptions
\ref{assume:alpha}, \ref{assume:alpha_upper} and equality
(\ref{error_model_d}) we can estimate derivative
$\dot{V}_{\hat{\thetavec}}$ as follows:
\begin{eqnarray}\label{parameric_deviation_derivative}
& &
\dot{V}_{\hat{\thetavec}}(\hat{\thetavec},\thetavec,t)\leq-(f(\bfx,\hat{\thetavec},t)-f(\bfx,\thetavec,t)+\varepsilon(t))(\hat{\thetavec}-\thetavec)^{T}\alphavec(\bfx,t)
- \frac{D}{4 D_1^2}\varepsilon^2(t) \nonumber
\\
& &
\leq-\frac{1}{D}(f(\bfx,\hat{\thetavec},t)-f(\bfx,\thetavec,t))^2+\frac{1}{D_1}|\varepsilon(t)||f(\bfx,\hat{\thetavec},t)-f(\bfx,\thetavec,t)|-
\frac{D}{4 D_1^2}\varepsilon^2(t)\\
\nonumber & & \leq -
\frac{1}{D}\left(|f(\bfx,\hat{\thetavec},t)-f(\bfx,\thetavec,t)|-\frac{D}{2
D_1} \varepsilon(t)\right)^2 \leq 0
\end{eqnarray}
It follows immediately from
(\ref{parameric_deviation_derivative}), (\ref{V_theta}) that
\begin{equation}\label{eq:parametric_norm}
\|\hat{\thetavec}(t)-\thetavec\|^{2}_{\Gamma^{-1}}\leq
\|\hat{\thetavec}(t_0)-\thetavec\|^{2}_{\Gamma^{-1}}+\frac{D}{2
D_1^2}\|\varepsilon(t)\|^{2}_{2,[t_0,\infty]}
\end{equation}
In particular, for $t\in[t_0,T^\ast]$ we can derive from
(\ref{V_theta}) that
\[
\|\hat{\thetavec}(t)-\thetavec\|^{2}_{\Gamma^{-1}}\leq
\|\hat{\thetavec}(t_0)-\thetavec\|^{2}_{\Gamma^{-1}}+\frac{D}{2
D_1^2}\|\varepsilon(t)\|^{2}_{2,[t_0,T^\ast]}
\]
Therefore $\hat{\thetavec}(t)\in L_\infty^2[t_0,T^\ast]$.
Furthermore
$|f(\bfx(t),\hat{\thetavec}(t),t)-f(\bfx(t),\thetavec,t)|-\frac{D}{2
D_1} \varepsilon(t)\in L_{2}^1 [t_0,T^\ast]$. In particular
\begin{equation}\label{eq:t1_ins1}
\begin{split}
\left\||f(\bfx(t),\hat{\thetavec}(t),t)-f(\bfx(t),\thetavec,t)|-\frac{D}{2
D_1} \varepsilon(t)\right\|_{2,[t_0,T^\ast]}^2\leq
\frac{D}{2}\|\thetavec-\hat{\thetavec}(t_0)\|^{2}_{\Gamma^{-1}}+\frac{D^2}{4
D_1^2}\|\varepsilon(t)\|_{2,[t_0,T^\ast]}^2
\end{split}
\end{equation}
Hence $f(\bfx(t),\hat{\thetavec}(t),t)-f(\bfx(t),\thetavec,t)\in
L_{2}^1 [t_0,T^\ast]$ as a sum of two functions from $L_{2}^1
[t_0,T^\ast]$. In order to estimate the upper bound of norm
$\|f(\bfx(t),\hat{\thetavec}(t),t)-f(\bfx(t),\thetavec,t)\|_{2,[t_0,T^\ast]}$
from (\ref{eq:t1_ins1}) we use the Minkowski inequality:
\[
\left\|f(\bfx(t),\hat{\thetavec}(t),t)-f(\bfx(t),\thetavec,t)|-\frac{D}{2
D_1} \varepsilon(t)\right\|_{2,[t_0,T^\ast]}\leq
\left(\frac{D}{2}\|\thetavec-\hat{\thetavec}(t_0)\|^{2}_{\Gamma^{-1}}\right)^{0.5}+
\frac{D}{2 D_1}\|\varepsilon(t)\|_{2,[t_0,T^\ast]}
\]
and then apply the triangle inequality to the functions from
$L_{2}^1 [t_0,T^\ast]$:
\begin{equation}\label{eq:t1_ins2}
\begin{split}
&
\|f(\bfx(t),\hat{\thetavec}(t),t)-f(\bfx(t),\thetavec,t)\|_{2,[t_0,T^\ast]}\leq
\left\|f(\bfx(t),\hat{\thetavec}(t),t)-f(\bfx(t),\thetavec,t)-\frac{D}{2
D_1}\varepsilon(t)\right\|_{2,[t_0,T^\ast]}+\\
& \frac{D}{2 D_1}\|\varepsilon(t)\|_{2,[t_0,T^\ast]}\leq
\left(\frac{D}{2}\|\thetavec-\hat{\thetavec}(t_0)\|^{2}_{\Gamma^{-1}}\right)^{0.5}
+ \frac{D}{D_1}\|\varepsilon(t)\|_{2,[t_0,T^\ast]}
\end{split}
\end{equation}
Therefore, property P1) is proven.

Let us prove property P2). In order to do this we have to check
first if the solutions of the closed loop system are defined for
all $t\in\Real_+$, i.e. they do not reach infinity in finite time.
We prove this by a contradiction argument. Indeed, let there
exists time instant $t_s$ such that $\|\bfx(t_s)\|=\infty$. It
follows from P1), however, that
$f(\bfx(t),\hat{\thetavec}(t),t)-f(\bfx(t),\thetavec,t)\in L_{2}^1
[t_0,t_s]$. Furthermore, according to (\ref{eq:t1_ins2}) the norm
$\|f(\bfx(t),\hat{\thetavec}(t),t)-f(\bfx(t),\thetavec,t)\|_{2,[t_0,t_s]}$
can be bounded from above by a continuous function of $\thetavec,
\ \hat{\thetavec}(t_0)$, $\Gamma$,  and
$\|\varepsilon(t)\|_{2,[t_0,\infty]}$. Let us denote this bound by
symbol $D_f$. Notice that $D_f$ does not depend on $t_s$. Consider
system (\ref{error_model_d}) for $t\in[t_0,t_s]$:
\[
\dpsi=f(\bfx,\thetavec,t)-f(\bfx,\hat{\thetavec},t)-\varphi(\psi,\omegavec,t)+\varepsilon(t)
\]
Given that both
$f(\bfx(t),\thetavec,t)-f(\bfx(t),\hat{\thetavec}(t),t),
\varepsilon(t) \in L_{2}^1 [t_0,t_s]$ and taking into account
Assumption \ref{assume:gain}, we automatically obtain that
$\psi(\bfx(t),t)\in L_\infty^{1}[t_0,t_s]$. In particular, using
the triangle inequality and the fact that function
$\gamma_{\infty,2}\left(\psi(\bfx_0,t_0),\omegavec,M\right)$ in
Assumption \ref{assume:gain} is non-decreasing in $M$, we can
estimate the norm $\|\psi(\bfx(t),t)\|_{\infty,[t_0,t_s]}$ as
follows:
\begin{equation}\label{eq:bound_psi}
\|\psi(\bfx(t),t)\|_{\infty,[t_0,t_s]}\leq
\gamma_{\infty,2}\left(\psi(\bfx_0,t_0),\omegavec,D_f+\|\varepsilon(t)\|^2_{2,[t_0,\infty]}\right)
\end{equation}
According to Assumption \ref{assume:psi} the following inequality
holds:
\begin{equation}\label{eq:bound_x}
\|\bfx(t)\|_{\infty,[t_0,t_s]}\leq\tilde{\gamma}\left(\bfx_0,\thetavec,\gamma_{\infty,2}\left(\psi(\bfx_0,t_0),\omegavec,D_f+\|\varepsilon(t)\|^2_{2,[t_0,\infty]}\right)\right)
\end{equation}
Given that a superposition of locally bounded functions is locally
bounded, we can conclude that $\|\bfx(t)\|_{\infty[t_0,t_s]}$ is
bounded.  This, however, contradicts to the previous claim that
$\|\bfx(t_s)\|=\infty$. Taking into account inequality
(\ref{eq:parametric_norm}) we can derive that both
$\hat{\thetavec}(\bfx(t),t)$ and $\hat{\thetavec}_I(t)$ are
bounded for every $t\in\Real_+$. Moreover, according to
(\ref{eq:bound_psi}), (\ref{eq:bound_x}),
(\ref{eq:parametric_norm}) these bounds are (locally bounded)
functions of initial conditions and parameters. Therefore,
$\bfx(t)\in L^n_\infty[t_0,\infty]$,
$\hat{\thetavec}(\bfx(t),t)\in L^d_\infty [t_0,\infty]$.
Inequality (\ref{eq:psi_gain}) follows immediately from
(\ref{eq:t1_ins2}), (\ref{eq:gain_psi_L2}), and the triangle
inequality.  Property P2) is proven.

Let us show that P3) holds. It is assumed that system
(\ref{eq:target_dynamics}) has $L_{2}^1 [t_0,\infty]\mapsto
L_{p}^1 [t_0,\infty]$, $p>1$ gain. In addition, we have just shown
that $f(\bfx(t),\thetavec,t)-f(\bfx(t),\hat{\thetavec}(t),t),
\varepsilon(t) \in L_{2} [t_0,\infty]$. Hence, taking into account
equation (\ref{error_model_d}) we conclude that
$\psi(\bfx(t),t)\in L_{p}^1 [t_0,\infty]$, $p>1$. On the other
hand, given that $f(\bfx,\hat{\thetavec},t)$,
$\varphi(\psi,\omegavec,t)$ are locally  bounded with respect to
their first two arguments uniformly in $t$, and that  $\bfx(t)\in
L_{\infty}^n[t_0,\infty]$,$\psi(\bfx(t),t)\in
L_\infty^1[t_0,\infty]$, $\hat{\thetavec}(t)\in
L_\infty^d[t_0,\infty]$, $\thetavec\in\Omega_\theta$, signal
$\varphi(\psi(\bfx(t),t),\omegavec,t)+f(\bfx(t),\thetavec,t)-f(\bfx(t),\hat{\thetavec}(t),t)$
is bounded. Then $\varepsilon(t)\in L_\infty^1[t_0,\infty]$
implies that $\dpsi$ is bounded, and P3) is guaranteed by
Barbalat's lemma.

To complete the proof of the theorem  (property P4) consider the
time-derivative of function $f(\bfx,\hat{\thetavec},t)$:
\[
\frac{d}{dt}f(\bfx,\hat{\thetavec},t)=L_{\bff(\bfx,\thetavec)+\bfg(\bfx)u(\bfx,\hat{\thetavec},t)}f(\bfx,\hat{\thetavec},t)+
\frac{\pd f(\bfx,\hat{\thetavec},t)}{\pd \hat{\thetavec}}\Gamma
(\varphi(\psi,\omegavec,t)+\dpsi)\alphavec(\bfx,t)+\frac{\pd
f(\bfx,\hat{\thetavec},t)}{\pd t}
\]
Taking into account that $\bff(\bfx,\thetavec)$, $\bfg(\bfx)$,
function $f(\bfx,\thetavec,t)$ is continuously differentiable in
$\bfx$, $\thetavec$; derivative $ \pd {f(\bfx,\thetavec,t)}/{\pd
t}$ is locally bounded with respect to $\bfx$, $\thetavec$
uniformly in $t$; functions $\alphavec(\bfx,t)$, $\pd
\psi(\bfx,t)/\pd t$ are locally bounded with respect to $\bfx$
uniformly in $t$, then $d/dt
(f(\bfx,\thetavec,t)-f(\bfx,\hat\thetavec,t))$ is bounded. Then
given that
$f(\bfx(t),\thetavec,t)-f(\bfx(t),\hat{\thetavec}(t),t)\in L_{2}^1
[t_0,\infty]$ by applying Barbalat's lemma we conclude that
$f(\bfx,\thetavec,\tau)-f(\bfx,\hat\thetavec,\tau)\rightarrow 0$
as $t\rightarrow\infty$. {\it The theorem is proven.}

\vskip 5mm

{\it Proof of Corollary \ref{cor:disturbance}.} Let
$\varepsilon(t)=0$. Then choosing the function
$V_{\hat{\thetavec}(\hat{\thetavec,\thetavec,t})}$ as in
(\ref{V_theta}), using (\ref{eq:dV_full_alg}), and invoking
Assumption \ref{assume:alpha}, we obtain that
\begin{equation}\label{dV_theta_e_zero}
\dot{V}_{\hat{\thetavec}(\hat{\thetavec},\thetavec,t)}\leq-(f(\bfx,\thetavec,t)-f(\bfx,\hat{\thetavec},t))\alphavec(\bfx,t)^{T}(\thetavec-\hat{\thetavec})\leq
-\frac{1}{D}(f(\bfx,\thetavec,t)-f(\bfx,\hat{\thetavec},t))^2
\end{equation}
Equality (\ref{dV_theta_e_zero}) and the fact that
$\varepsilon(t)=0$ in (\ref{V_theta}) imply that the norm
$\|\hat{\thetavec}-\thetavec\|^2_{\Gamma^{-1}}$ is non-increasing.
Furthermore, (\ref{dV_theta_e_zero}) implies that
\begin{equation}\label{eq:diff_f_2_e_zero}
\|f(\bfx(t),\thetavec,t)-f(\bfx(t),\hat{\thetavec}(t),t)\|_{2,
[t_0,T^\ast]}\leq
\left(\frac{D}{2}\|\hat{\thetavec}(t_0)-\thetavec\|^2_{\Gamma^{-1}}\right)^{0.5}
\end{equation}
This proves property P1). Taking into account
(\ref{eq:diff_f_2_e_zero}) and given that Assumptions
\ref{assume:psi}, \ref{assume:gain} are satisfied we can conclude
that $\bfx(t)\in L_{\infty}^n[t_0,\infty]$, $\psi(\bfx(t),t)\in
L_\infty^1[t_0,\infty]$, and the following estimate holds:
\begin{equation}\label{eq:psi_gain_e_zero}
\|\psi(\bfx(t),t)\|_{\infty,[t_0,\infty]}\leq
\gamma_{\infty,2}\left(\psi(\bfx_0,t_0),\omegavec,\left(\frac{D}{2}\|\hat{\thetavec}(t_0)-\thetavec\|^2_{\Gamma^{-1}}\right)^{0.5}
\right)
\end{equation}
Hence P2) is also proven. Properties P3),P4) follow by the same
arguments as in the proof of Theorem \ref{stability_theorem}.
Therefore, P5) is proven.
%
%
%
%In order to show that P6) holds we have to make sure that
%$f(\bfx(t),\thetavec,t)-f(\bfx(t),\hat{\thetavec}(t),t)\in
%L_2^{1}[t_0,\infty]$  for all $\alphavec(\bfx,t): \
%\|\alphavec(\bfx,t)\|\leq D_\alpha$ without invoking Assumption
%\ref{assume:alpha_upper}. To this objective consider the following
%positive-definite function
%$V_{0,\hat{\thetavec}}(\hat{\thetavec},\thetavec)=0.5\|\hat{\thetavec}-\thetavec\|^2_{\Gamma^{-1}}$.
%Its time-derivative can be estimated in the same manner as in
%(\ref{parameric_deviation_derivative}):
%\[
%\dot{V}_{0,\hat{\thetavec}}\leq -
%\]
{\it The corollary is proven.}

\vskip 5mm

{\it Proof of Corollary \ref{cor:gains}.} Let us show that P6)
holds. Without the loss of generality assume that solutions of the
system exist over the following time interval $[t_0,T^\ast]$.
According to Theorem \ref{stability_theorem}, property P1), the
norm $\|\thetavec-\hat{\thetavec}(t)\|$ is bounded from above by a
function of initial conditions $\hat{\thetavec}(t_0)$, parameters
$\Gamma$, $D$, $D_1$, and $\|\varepsilon(t)\|_{2,[t_0,\infty]}$.
Let us denote this bound by symbol $B_\theta$. Notice that
$B_\theta$ does not depend on $T^\ast$. On the other hand,
according to Hypothesis H\ref{hyp:globally_Lipschitz_uniform} the
following estimate holds:
\[
\exists \ D_\theta>0:  \
|f(\bfx,\thetavec,t)-f(\bfx,\hat{\thetavec},t)|\leq
D_\theta\|\thetavec-\hat{\thetavec}\|
\]
Hence $|f(\bfx(t),\thetavec,t)-f(\bfx(t),\hat{\thetavec}(t),t)|\in
L_\infty^1[t_0,T^\ast]$ and moreover
\begin{equation}\label{eq:cor_df_infty}
\|f(\bfx(t),\thetavec,t)-f(\bfx,\hat{\thetavec}(t),t)\|_{\infty,[t_0,T^\ast]}\leq
D_\theta B_\theta
\end{equation}
Consider the following signal
$\mu(t)=f(\bfx(t),\thetavec,t)-f(\bfx(t),\hat{\thetavec}(t),t)+\varepsilon(t)$.
Signal $\mu(t)\in L_{2}^1[t_0,T^\ast]\cap L_\infty^1[t_0,T^\ast]$,
and let $M_\infty$,  $M_2\in\Real_+$ be the bounds for its
$L_\infty^1[t_0,T^\ast]$, $L_2^1[t_0,T^\ast]$ norms respectively.
According to (\ref{eq:t1_ins2}), (\ref{eq:cor_df_infty}), these
bounds can be estimated as follows: $M_\infty=D_\theta
B_\theta+\|\varepsilon(t)\|_{\infty,[t_0,\infty]}$,
$M_2=\left(\frac{D}{2}\|\thetavec-\hat{\thetavec}(t_0)\|^{2}_{\Gamma^{-1}}\right)^{0.5}
+
\left(\frac{D}{D_1}+1\right)\|\varepsilon(t)\|_{2,[t_0,\infty]}$.
Therefore, given $p\geq 2$ one can derive that
\[
\int_{t_0}^{\infty}\mu^p(\tau)d\tau=\int_{t_0}^{\infty}\mu^{p-2}(\tau)\mu^{2}(\tau)d\tau\leq
M^{p-2}_\infty M_2^2
\]
Hence, $\mu(t)\in L_{p}^1[t_0,T^\ast]$ and its
$L_p[t_0,T^\ast]$-norm is bounded from above by $M^{p-2}_\infty
M_2^2$, where the bounds $M_\infty$, $M_2$ both do not depend on
$T^\ast$. According to the corollary formulation, system
(\ref{eq:target_dynamics}) has $L_{p}^1[t_0,\infty]\rightarrow
L_\infty^1[t_0,\infty]$ gain and therefore $\psi(\bfx(t),t)\in
L_\infty[t_0,T^\ast]$. Then applying the same argument as in the
proof of property P2) of Theorem \ref{stability_theorem} and using
Assumption \ref{assume:psi} we can immediately obtain that
$\bfx(t)\in L_\infty^n[t_0,\infty]$, $\psi(\bfx(t),t)\in
L_\infty^1[t_0,\infty]$, and $\hat{\thetavec}(t)\in
L_{\infty}^d[t_0,\infty]$. Thus property P6) is proven. Property
P7) can now be proven in the same way as property P3) in Theorem
\ref{stability_theorem}. {\it The corollary is proven.}

\vskip 5mm

{\it Proof of Theorem \ref{convergence_theorem}.} According to the
theorem formulation, Assumptions
\ref{assume:psi},\ref{assume:gain}, \ref{assume:alpha},
\ref{assume:explicit_realizability} hold. Hence, applying
Corollary \ref{cor:disturbance} we can conclude that
$\hat{\thetavec}(t)\in L_\infty^d[t_0,\infty]$ and $\bfx(t)\in
L_\infty^n[t_0,\infty]$.

Let us show that limiting relation (\ref{eq:convergence}) holds in
case alternative 1) is satisfied. To this purpose consider
derivative $\dot{\hat{\thetavec}}$:
\begin{equation}\label{eq:equivalent_alg_1}
\dot{\hat{\thetavec}}=\Gamma(\dpsi+\varphi(\psi))\alphavec(\bfx,t)=\Gamma
(f(\bfx,\thetavec,t)-f(\bfx,\hat{\thetavec},t))\alphavec(\bfx,t).
\end{equation}
Given that $\hat{\thetavec}(t)\in L_\infty^d[t_0,\infty]$ and
$\bfx(t)\in L_\infty^n[t_0,\infty]$, and that Hypothesis
H\ref{hyp:local_linear_bound} holds, the function
$f(\bfx,\thetavec,t)$ satisfies the following inequality for some
$D$, $D_1, \in \Real_+$:
\begin{eqnarray}
& &
D_1|\alphavec(\bfx,t)^{T}(\hat{\thetavec}-\thetavec)|\leq|f(\bfx,\thetavec,t)-f(\bfx,\hat{\thetavec},t))|\leq
D|\alphavec(\bfx,t)^{T}(\hat{\thetavec}-\thetavec)|;\nonumber\\
& &
\alphavec(\bfx,t)^{T}(\hat{\thetavec}-\thetavec)(f(\bfx,\hat\thetavec,t)-f(\bfx,\thetavec,t))\geq0
\nonumber
\end{eqnarray}
Therefore, there exists function
$\kappa:\Real_+\rightarrow\Real_+$, $D_1\leq\kappa^2(t)\leq D$
such that
\begin{equation}\label{eq:parameters_eq}
\dot{\hat{\thetavec}}=-\kappa^2(t)\Gamma
\alphavec(\bfx,t)^{T}(\hat{\thetavec}-\thetavec)\alphavec(\bfx,t)=-\kappa^2(t)\Gamma
\alphavec(\bfx,t)\alphavec(\bfx,t)^{T} (\hat{\thetavec}-\thetavec)
\end{equation}
Notice that matrix $\Gamma$ is a positive definite and symmetric
matrix. It, therefore, can be factorized as follows
$\Gamma=\Gamma_0 \Gamma_0^T$:  where $\Gamma_0$ is nonsingular
$n\times n$ real matrix. Let us define
$\tilde{\thetavec}=\Gamma_0^{-1}(\hat{\thetavec}-\thetavec)$. In
these new coordinates, equation (\ref{eq:parameters_eq}) will have
the following form:
\begin{equation}\label{eq:parameters_eq_2}
\dot{\tilde{\thetavec}}=-\kappa(t)^2\Gamma_0^{-1}\Gamma_0\Gamma_0^T
\alphavec(\bfx,t)\alphavec(\bfx,t)^T(\hat{\thetavec}-\thetavec)=-\kappa^2(t)\Gamma_0^T\alphavec(\bfx,t)\alphavec(\bfx,t)^T\Gamma_0\tilde\thetavec
\end{equation}
Denoting $\kappa(t)\Gamma_0^T\alphavec(\bfx,t)=\phivec(\bfx,t)$ we
can rewrite equation (\ref{eq:parameters_eq_2}) as follows:
\begin{equation}\label{eq:parameters_gradient_flow}
\dot{\tilde{\thetavec}}=-\phivec(\bfx,t)\phivec(\bfx,t)^T\tilde{\thetavec}
\end{equation}
where function
$\phivec(\bfx,t):\Real^n\times\Real_+\rightarrow\Real^d$ satisfies
equality:
\begin{eqnarray}\label{eq:parameters_eq_PE}
& &
\etavec^T\int_t^{t+L}\phivec(\bfx(\tau),\tau)\phivec(\bfx(\tau),\tau)d\tau\etavec=\etavec^T\Gamma_0^T\left(\int_{t}^{t+L}\kappa^2(\tau)\alphavec(\bfx(\tau),\tau)\alphavec(\bfx(\tau),\tau)^{T}d\tau\right)\Gamma^0\etavec
\end{eqnarray}
for all $\etavec\in\Real^d$. Taking into account that function
$\alphavec(\bfx(t),t)$ is persistently exciting,
$\Gamma=\Gamma_0^T\Gamma_0$, and that $\kappa^2(t)\geq D_1$ we can
obtain the following bound for quadratic form
(\ref{eq:parameters_eq_PE}):
\begin{eqnarray}\label{eq:parameters_eq_new_PE}
& &
\etavec^T\int_t^{t+L}\phivec(\bfx(\tau),\tau)\phivec(\bfx(\tau),\tau)d\tau\etavec\geq
\delta D_1 \|\Gamma_0\etavec\|^2=\delta D_1 \etavec^T \Gamma
\etavec \geq \delta D_1 \lambda_{\min}(\Gamma)
\|\etavec\|^2=\delta_\phi \|\etavec\|^2
\end{eqnarray}
Hence, function $\phivec(\bfx(t),t)$ is also persistently
exciting. Notice also that $\|\phi(\bfx,t)\|$ is bounded from
above:
\begin{equation}
\begin{split}
\|\phivec(\bfx,t)\|&=\|\kappa(t)\alphavec(\bfx,t)\Gamma_0\|\leq{\lambda_{\max}(\Gamma_0)}\|\kappa(t)\alphavec(\bfx,t)\|\leq
{\lambda_{\max}(\Gamma_0) D} \alpha_\infty  \nonumber \\
\alpha_\infty & =\sup_{\|\bfx\|\leq
\|\bfx(t)\|_{\infty,[t_0,\infty]}, \ t\geq t_0}
\|\alphavec(\bfx,t)\|
\end{split}
\end{equation}
In order to show that
$\lim_{t\rightarrow\infty}\tilde{\thetavec}(t)=0$ exponentially
fast we invoke the following useful lemma from \cite{Lorea_2002}
(Lemma 5, page 18):
\begin{lem}\label{lem:Lorea_gradient_flow} Let system
(\ref{eq:parameters_gradient_flow}) be given, condition
(\ref{eq:parameters_eq_new_PE}) holds (uniformly), and
$\phivec(\bfx,t)$ in (\ref{eq:parameters_gradient_flow}) be
bounded $\|(\phivec(\bfx,t))\|\leq \phi_M$. Then system
(\ref{eq:parameters_gradient_flow}) is (uniformly) exponentially
stable and, furthermore:
\begin{equation}\label{eq:parameters_rate_0}
\|\tilde{\thetavec}(t)\|\leq e^{- K t}\|\tilde{\thetavec}(t_0)\|,
\ K=\frac{\delta_\phi}{L(1+L\phi_M^2)^2}
\end{equation}
\end{lem}
According to Lemma \ref{lem:Lorea_gradient_flow} solutions of
system (\ref{eq:parameters_gradient_flow}) converge to the origin
exponentially fast with a rate of convergence defined by
(\ref{eq:parameters_rate_0}), where
\begin{equation}\label{eq:parameters_rate_1}
\delta_\phi=\lambda_{\min}(\Gamma)D_1\delta, \
\phi_M^2=\lambda_{\max}(\Gamma)D^2\alpha_\infty^2
\end{equation}
Taking into account equation (\ref{eq:parameters_rate_0}),
(\ref{eq:parameters_rate_1}) and observing that
$(1+L\phi_M^2)^2\leq 2(1+\phi_M^4 L^2)$,
$\lambda_{\max}(\Gamma_0)^2=\lambda_{\max}(\Gamma)$ we can
estimate $\|\tilde{\thetavec}(t)\|$ as follows:
\begin{equation}\label{eq:parameters_rate_2}
\|\tilde{\thetavec}(t)\|\leq e^{- K_1
t}\|\tilde{\thetavec}(t_0)\|, \ K_1=\frac{\delta
D_1\lambda_{\min}(\Gamma)}{2L(1+\lambda_{\max}^2(\Gamma)L^2 D^2
\alpha_\infty^4)}
\end{equation}
Given that $\|\Gamma_0
\tilde{\thetavec}(t)\|=\|(\hat{\thetavec}(t)-\thetavec)\|$, and
using (\ref{eq:parameters_rate_2}) we derive the following bounds
for $\|(\hat{\thetavec}(t)-\thetavec)\|$:
\[
\|(\hat{\thetavec}(t)-\thetavec)\|\leq
\|\Gamma_0\|\|\tilde{\thetavec}(t)\|\leq
{\lambda_{\max}(\Gamma_0)}e^{- K_1
t}\|\Gamma_0^{-1}(\hat{\thetavec}(t_0)-\thetavec)\|\leq
\left(\frac{\lambda_{\max}(\Gamma_0)}{\lambda_{\min}(\Gamma_0)}\right)e^{-K_1
t}\|\hat{\thetavec}(t_0)-\thetavec\|
\]
This proves alternative 1) of the theorem.

Let us prove alternative 2). It follows immediately from Corollary
\ref{cor:disturbance} of Theorem \ref{stability_theorem} that
\begin{equation}\label{eq:function_difference}
\lim_{t\rightarrow\infty}f(\bfx(t),\thetavec,t)-f(\bfx(t),\hat{\thetavec}(t),t)=0
\end{equation}
Furthermore, given that
$\dot{\hat{\thetavec}}=\Gamma(f(\bfx(t),\thetavec,t)-f(\bfx(t),\hat{\thetavec}(t),t))\alphavec(\bfx,t)$,
$\bfx(t)\in L_\infty^n[t_0,\infty]$ and $\alphavec(\bfx,t)$ is
locally bounded in $\bfx$ uniformly in $t$, we can conclude that
$\dot{\hat{\thetavec}}\rightarrow 0$ as $t\rightarrow\infty$. Let
us divide the $\Real_+$ into the following union of subintervals:
\[
\Real_+=\bigcup_{i=1}^\infty \Delta_i, \ \Delta_i=[t_i,t_i+T], \
t_0=0, \ t_{i+1}=t_{i}+T, \ i\in\Natural
\]
The fact that $\dot{\hat{\thetavec}}\rightarrow 0$ as
$t\rightarrow\infty$ ensures that
\begin{equation}\label{eq:parameter_convergence_interval}
\lim_{i\rightarrow\infty}\|\hat{\thetavec}(s_i)-\hat{\thetavec}(\tau_i)\|=0,
\ \forall \  s_i,\tau_i\in\Delta_i
\end{equation}
In order to show this let us integrate equation
(\ref{eq:equivalent_alg_1})
\begin{equation}\label{eq:parameter_difference_estimate_1}
\|\hat{\thetavec}(s_i)-\hat{\thetavec}(\tau_i)\|=\|\Gamma
\int_{s_i}^{\tau_i}(f(\bfx(\tau),{\thetavec},\tau)-f(\bfx(\tau),\hat{\thetavec}(\tau),\tau))\alpha(\bfx(\tau),\tau)d\tau\|
\end{equation}
Applying the Cauchy-Schwartz inequality to
(\ref{eq:parameter_difference_estimate_1}) and subsequently using
the mean value theorem we can obtain the following estimate:
\begin{eqnarray}\label{eq:parameter_difference_estimate_2}
\|\hat{\thetavec}(s_i)-\hat{\thetavec}(\tau_i)\|&\leq&\int_{t_i}^{t_i+T}\|\Gamma\|\cdot |f(\bfx(\tau),{\thetavec},\tau)-f(\bfx(\tau),\hat{\thetavec}(\tau),\tau)|\cdot\|\alpha(\bfx(\tau),\tau)\|d\tau\nonumber\\
&= & \|\Gamma\| \cdot T \cdot
 |f(\bfx(\tau_i'),\thetavec,\tau_i')-f(\bfx(\tau_i'),\hat{\thetavec}(\tau_i'),\tau_i')|\cdot \|\alphavec(\bfx(\tau_i'),\tau_i')\|,
\ \tau_i'\in\Delta_i
\end{eqnarray}
Given that limiting relation (\ref{eq:function_difference}) holds,
$\bfx(t)\in L_\infty^n[t_0,\infty]$, and $\alphavec(\bfx,t)$ is
locally bounded uniformly in $t$ we can conclude from
(\ref{eq:parameter_difference_estimate_2}) that limiting relation
(\ref{eq:parameter_convergence_interval}) holds.

Let us choose a sequence of points from $\Real_+$:
$\{\tau_i\}_{i=1}^\infty$ such that $\tau_i\in\Delta_i$,
$i\in\Natural$. As follows from the nonlinear persistent
excitation condition (inequality (\ref{eq:NLPE_2})), for every
$\hat{\thetavec}(\tau_i)$, $\tau_i\in\Delta_i$ there exists a
point $t'_i\in\Delta_i$ such that the following inequality holds
\begin{equation}\label{eq:parameters_nonlinear_PE}
\|f(\bfx(t'_i),\thetavec,t'_i)-f(\bfx(t'_i),\hat{\thetavec}(\tau_i),t'_i)\|\geq
\varrho(\|\thetavec-\hat{\thetavec}(\tau_i)\|)\geq 0
\end{equation}
Let us consider the following differences:
\[
f(\bfx(t'_i),\hat{\thetavec}(\tau_i),t'_i)-f(\bfx(t'_i),\hat{\thetavec}(t'_i),t'_i),
\ \tau_i,t'_i\in\Delta_i
\]
It follows immediately from H\ref{hyp:locally_bound_uniform_f},
H\ref{hyp:locally_bound_uniform_df}, and
(\ref{eq:parameter_convergence_interval}) that
\begin{equation}\label{eq:function_difference_residual}
\lim_{i\rightarrow\infty}
f(\bfx(t'_i),\hat{\thetavec}(\tau_i),t'_i)-f(\bfx(t'_i),\hat{\thetavec}(t'_i),t'_i)=0,
\  \ \tau_i,t'_i\in\Delta_i
\end{equation}
Taking into account (\ref{eq:function_difference_residual}) and
(\ref{eq:function_difference}) we can derive that
\begin{eqnarray}\label{eq:nonlinear_pe_convergence}
& & \lim_{i\rightarrow\infty}
f(\bfx(t'_i),\thetavec,t'_i)-f(\bfx(t'_i),\hat{\thetavec}(\tau_i),t'_i)=\lim_{i\rightarrow\infty}(f(\bfx(t'_i),\thetavec,t'_i)-f(\bfx(t'_i),\hat{\thetavec}(t'_i),t'_i))+\nonumber\\
& &
\lim_{i\rightarrow\infty}f(\bfx(t'_i),\hat{\thetavec}(t'_i),t'_i)-f(\bfx(t'_i),\hat{\thetavec}(\tau_i),t'_i)=0
\end{eqnarray}
According to (\ref{eq:nonlinear_pe_convergence}) and
(\ref{eq:parameters_nonlinear_PE}), sequence
$\{\varrho(\|\thetavec-\hat{\thetavec}(\tau_i)\|)\}_{i=1}^\infty$
is bounded from above and below by two sequences converging to
zero. Hence,
$\lim_{i\rightarrow\infty}\varrho(\|\thetavec-\hat{\thetavec}(\tau_i)\|)=0$.
Notice that $\varrho(\cdot)\in \mathcal{K}\cap C^0$ which implies
that
\begin{equation}\label{eq:parameter_convergence}
\lim_{i\rightarrow\infty}\|\thetavec-\hat{\thetavec}(\tau_i)\|=0
\end{equation}
In order to show that
$\lim_{t\rightarrow\infty}({\thetavec}-\hat{\thetavec}(t))=0$
notice that
\[
\|\thetavec-\hat{\thetavec}(t)\|\leq
\|\thetavec-\hat{\thetavec}(s_i)\|, \
s_i=\arg\max_{s\in\Delta_i}\|\thetavec-\hat{\thetavec}(s)\| \ \
\forall \ t\in\Delta_i
\]
Hence, applying the triangle inequality
$\|\thetavec-\hat{\thetavec}(s_i)\|\leq
\|\thetavec-\hat{\thetavec}(\tau_i)\|+\|\hat{\thetavec}(\tau_i)-
\hat{\thetavec}(s_i)\|$ and using equations
(\ref{eq:parameter_convergence_interval}),
(\ref{eq:parameter_convergence}) we can conclude that
$\|\thetavec-\hat{\thetavec}(t)\|$ is bounded from above and below
by two functions converging to zero. Hence,
$\|\thetavec-\hat{\thetavec}(t)\|\rightarrow 0$ as $t\rightarrow
\infty$ and limiting relation (\ref{eq:convergence}) holds.  {\it
The theorem is proven.}

\bibliographystyle{plain}
\bibliography{identification_fifo_ver4}

\end{document}